\documentclass{article}
\usepackage{latexsym,amsmath,amscd,amssymb,framed,graphics}
\usepackage{enumerate}
\usepackage{graphicx}
\usepackage{amsthm}
\usepackage{epstopdf}

\usepackage[margin=2cm]{geometry}
\usepackage{palatino}
\usepackage[colorlinks]{hyperref}

\def\RR{\mathbb R}
\def\CC{\mathbb C}

\newcommand \al{\alpha}
\newcommand\be{\beta}

\newcommand\om{\omega}

\newcommand\ie{i.e.\ }
\newcommand\F{\mathcal{F}}

\newcommand\g{\mathfrak g}
\newcommand\x{\times}

\newcommand\pa{\partial}
\renewcommand\o{\circ}
\newcommand\on{\operatorname}
\newcommand\su{\mathfrak {su}}
\newcommand\SU{\on{SU}}
\renewcommand\u{\mathfrak {u}}
\newcommand\U{\on{U}}
\newcommand\rank{\on{rank}}
\newcommand\Ad{\on{Ad}}
\renewcommand\v{\mathbf{v}}
\newcommand\w{\mathbf{w}}
\renewcommand\a{\mathbf{a}}
\renewcommand\b{\mathbf{b}}

\begin{document}

\newtheorem{theorem}{Theorem}[section]
\newtheorem{definition}[theorem]{Definition}
\newtheorem{lemma}[theorem]{Lemma}
\newtheorem{remark}[theorem]{Remark}
\newtheorem{proposition}[theorem]{Proposition}
\newtheorem{corollary}[theorem]{Corollary}
\newtheorem{example}[theorem]{Example}

\title{The n:m resonance dual pair} 
\author{Darryl D. Holm$^{1}$ and Cornelia Vizman$^{2}$ }
\addtocounter{footnote}{1} 
\footnotetext{Department of Mathematics, Imperial College London. London SW7 2AZ, UK. Partially supported by Royal Society of London Wolfson Award.
\texttt{d.holm@imperial.ac.uk}
\addtocounter{footnote}{1} }
\footnotetext{Department of Mathematics,
West University of Timi\c soara, 300223 Timi\c soara, Romania.
\texttt{vizman@math.uvt.ro}
\addtocounter{footnote}{1} }
\date{In honor of Tudor Ratiu's sixtieth birthday. }

\maketitle
\makeatother

\maketitle

\noindent \textbf{AMS Classification:} 53D17, 53D20

\noindent \textbf{Keywords:} dual pairs, Poisson brackets

\begin{abstract}
\end{abstract}


\section{Introduction}\label{sec1}

The dual pairs for 1:1 and 1:-1 resonance are presented in \cite{Ma},
reformulating results from \cite{CR} and \cite{Iwai}.
The dual pair for 1:1 resonance is the pair 
of  momentum maps associated to the commuting Hamiltonian actions 
of the Lie groups $S^1$ and ${\SU(2)}$ on  $\CC^2$ endowed with  the opposite $\om$ of
the canonical symplectic form:
\[
{\RR\stackrel{R}{\longleftarrow}(\CC^2,\om)
\stackrel{J}{\longrightarrow} \su(2)^*}.
\] 
The momentum map $J$ maps the fibers of $R$, which are 3-spheres, 
into 2-spheres, coadjoint orbits of $\SU(2)$.
The restriction of $J$ to these 3-spheres is a Hopf fibration.

A similar construction works for the $S^1$ and $\SU(1,1)$ actions on $\CC^2$ endowed 
with the symplectic form $\om_-=-dx_1\wedge dy_1+dx_2\wedge dy_2$, 
thus obtaining the 1:-1 resonance dual pair: 
\[
{\RR\stackrel{R_-}{\longleftarrow}(\CC^2,\om_-)
\stackrel{J_-}{\longrightarrow} \su(1,1)^*}.
\] 
The momentum map $J_-$ maps the fibers of $R_-$, which are 3-hyperboloids, 
into 2-hyperboloids, coadjoint orbits of $\SU(1,1)$.
The restriction of $J_-$ to these 3-hyperboloids is a hyperbolic Hopf fibration.

In this paper we build dual pairs of Poisson maps
\[
\RR\stackrel{R_\pm}{\longleftarrow}(D,\om_\pm)
\stackrel{\Pi_\pm}{\longrightarrow} B
\] 
associated to $n:m$ resonance, as well as to $n:-m$ resonance.
Except for the above mentioned cases $1:\pm1$, these are not pairs of momentum maps.
Here $D$ is an open subset of $\CC^2$ with the above mentioned symplectic forms 
$\om_\pm$, and $B$ an open subset of $\RR^3$.
The Poisson structure on $B$, which depends on the natural numbers $n$ and $m$,  is not Lie-Poisson.
Instead, its symplectic leaves are the Kummer shapes: 
bounded surfaces for $n:m$ resonance, 
and unbounded surfaces for $n:-m$ resonance \cite{Ku1986}.

Under some extra hypothesis, to each integrable system in the non-commutative sense (also called superintegrable system) one can associate a dual pair
whose right leg is the map defined by the independent first integrals \cite{Fasso} \cite{OrRa04}.
Beside the rigid body and the Kepler system,
the two uncoupled oscillators in $m:n$ resonance comprise a well known example of superintegrable system.
The dual pairs we present in this article are of this type.

\paragraph{Acknowledgements.}
We are grateful to Andreas Kriegl for very helpful suggestions regarding Lemma \ref{hard} and Lemma \ref{hardbis},
and to Francesco Fasso for turning our attention to dual pairs in the context of superintegrable systems. We also acknowledge partial support by the Royal Society of London's Wolfson Scheme and  hospitality at the Institute for Mathematical Sciences, Imperial College London. Finally, we are grateful to our late friend Jerry Marsden and we fondly remember the many wonderful discussions of geometric mechanics we had together with both Jerry and Tudor over the years.


\section{Dual pairs}\label{sec2}

Let $(M,\omega)$ be a symplectic manifold and $P_1, P_2$ be two Poisson manifolds. A pair of Poisson mappings
\begin{equation*}
P_1\stackrel{J_1}{\longleftarrow}(M,\om)\stackrel{J_2}{\longrightarrow} P_2
\end{equation*}
is called a \textit{dual pair} \cite{We83}
if $\ker TJ_1$ and $\ker TJ_2$ are symplectic orthogonal complements of one another. That is
\begin{equation}\label{finite_dimensional_dual_pair}
(\ker TJ_1)^\om=\ker TJ_2.
\end{equation}
A systematic treatment of dual pairs can be found in Chapter 11 of \cite{OrRa04}.
The infinite dimensional case is treated in \cite{GBV10}.

\begin{proposition}\label{propo}
Let $J_1$ and $J_2$ be momentum maps arising from the canonical actions of two connected
Lie groups $G_1$ and $G_2$ on a symplectic manifold $(M,\om)$. 
We assume that both momentum maps are equivariant, so they are Poisson maps with
respect to the (+) Lie-Poisson structure on the dual Lie algebras.
Moreover we assume that $J_1$ is $G_2$-invariant,
and the $G_2$ action is transitive on level sets of $J_1$.
Then the pair of momentum maps
\begin{equation}\label{dp}
\g_1^*\stackrel{J_1}{\longleftarrow}(M,\om)\stackrel{J_2}{\longrightarrow}\g_2^*
\end{equation}
is a dual pair.
\end{proposition}

\begin{proof}
The transitivity of the $G_2$ action on level sets of $J_1$
is written infinitesimally as $\left(\mathfrak{g}_2\right)_M=\ker T{\bf J}_1$.
The dual pair property is seen upon writing $(\ker TJ_1)^\om=((\g_2)_M)^\om=\ker TJ_2$.
\end{proof}

As a consequence we obtain that the actions of the Lie groups $G_1$ and $G_2$ on $M$ commute.

The dual pair is called \textit{full} if $J_1:M\rightarrow P_1$ and $J_2:M\rightarrow P_2$ are surjective submersions. 
A key result in the context of dual pairs is the symplectic leaf correspondence for full dual pairs with connected fibers. Namely, there is a bijective
correspondence between the symplectic leaves 
of $P_1$ and those 
of $P_2$ \cite{We83}:
$$\mathcal{L}_1\mapsto J_2(J_1^{-1}(\mathcal{L}_1))\text{ with inverse }
\mathcal{L}_2\mapsto J_1(J_2^{-1}(\mathcal{L}_2)).$$


\section{The $1:1$ resonance dual pair}\label{sec3}

Let $\langle\ ,\ \rangle$ be the canonical Hermitian inner product on $\CC^2$.
This means that $\langle\ ,\ \rangle=g(\ ,\ )+i\om(\ ,\ )$, 
with $g$ the euclidean metric on $\CC^2$ and $\om$ the opposite of the canonical symplectic form on $\CC^2$.  
The Lie group $\U(2)$ of unitary $2\x 2$ matrices, \ie complex matrices $g$ with the property
$\langle g\a,g\b\rangle=\langle\a,\b\rangle$ for all $\a,\b\in\CC^2$,
acts in a Hamiltonian way on $(\CC^2,\om)$ with momentum map
\begin{equation}\label{canon}
\bar J:\CC^2\to\u(2)^*,\quad \langle\bar J(\a),\xi\rangle_{\u(2)}=\frac i2\langle\a,\xi(\a)\rangle.
\end{equation}
This follows from the computation
\begin{align*}
d_\a\langle\bar J,\xi\rangle_{\u(2)}=\frac i2\langle\a,\xi(\ )\rangle+\frac i2\langle\ ,\xi(\a)\rangle
=\Im\langle\xi(\a),\ \rangle=\om(\xi_{\CC^2}(\a),\ ),
\end{align*}
where $\xi\in\u(2)$, the Lie algebra of skew Hermitian $2\x 2$ matrix, \ie $\langle\xi(\a),\b\rangle+\langle\a,\xi(\b)\rangle=0$
for all $\a,\b\in\CC^2$.
Another way to deduce this is from the general form of the momentum map for linear 
Hamiltonian actions on linear symplectic
spaces $\langle\bar J(\a),\xi\rangle_{\u(2)}=\frac 12\om(\xi(\a),\a)$, since $g(\xi(\a),\a)=0$
for all $\xi\in\u(2)$.

The Lie group $\U(2)$ is the direct product of its center,
which is isomorphic to the circle $S^1$,
and the special unitary group 
$$
\SU(2)=\left\{\left[\begin{array}{ccc} \al & \be \\ -\bar\be & \bar\al \end{array}\right]
\text{ with }|\al|^2+|\be|^2=1\right\}.
$$
 
The momentum map for the circle action is 
\begin{align*}
R:\CC^2\to\RR,\quad R(\a)&=\frac 12(|a_1|^2+|a_2|^2)=\frac12\langle\a,\a\rangle.
\end{align*}
To compute the momentum map for the $\SU(2)$ action we consider the linear isomorphism
\begin{align*}
v\in\RR^3&\mapsto 
\xi_v=\left[\begin{array}{ccc} iv_3 & iv_1+v_2 \\ iv_1-v_2 & -iv_3 \end{array}\right]\in\su(2)\subset\u(2).
\end{align*}
With this identification, using again the expression \eqref{canon} of $\bar J$,
the momentum map becomes
\begin{equation}\label{j_can}
J:\CC^2\to \su(2)^*=\RR^3,\quad J(\a)=\left({\rm Re}(a_1 \bar a_2) , -\,{\rm Im}(a_1 \bar a_2) ,
\frac 12|a_1|^2-\frac 12|a_2|^2 \right).
\end{equation}
This follows from the computation:
\begin{align*}
\langle J(\a),v\rangle_{\RR^3}&=\langle \bar J(\a),\xi_v\rangle_{\u(2)}
=\frac i2\langle\a,\xi_v(\a)\rangle=\frac 12 v_1(a_1\bar a_2+\bar a_1a_2)+\frac i2 v_2(a_1\bar a_2-\bar a_1 a_2)
+\frac 12 v_3(a_1\bar a_1-a_2\bar a_2).
\end{align*}

It is easy to see that the momentum map $\bar J$ is $\U(2)$-equivariant:
\[
\langle\bar J(g\cdot\a,\xi\rangle_{\u(2)}
=\frac12i\langle g\cdot\a,\xi(g\cdot\a)\rangle
=\frac12i\langle\a,(g^{-1}\xi g)(\a)\rangle
=\langle\bar J(\a),\Ad_{g^{-1}}\xi\rangle_{\u(2)}
=\langle\Ad^*_{g^{-1}}\bar J(\a),\xi\rangle_{\u(2)}.
\]
Therefore, ${\bar J}(g\cdot\mathbf{a}) = \Ad^*_{g^{-1}}{\bar J}(\mathbf{a})$
for all $g\in\U(2)$.
From the equivariance of $\bar J$ follows the $S^1$ equivariance of $R$ and the $\SU(2)$ equivariance of $J$.

We denote the components of the momentum map $J$ by $X,Y,Z$, so
\[
X(\a)-iY(\a)= a_1 \bar a_2
\quad\hbox{and}\quad
Z(\a)= \frac 12|a_1|^2-\frac 12|a_2|^2.
\]
They satisfy $X^2 + Y^2 + Z^2 = R^2$. In real coordinates we recognize the three first integrals of the integrable system of two uncoupled oscillators:
\begin{align*}
X (x_1,y_1,x_2,y_2) & =x_1x_2+y_1y_2\\
Y(x_1,y_1,x_2,y_2) & =x_1y_2-x_2y_1\\
Z(x_1,y_1,x_2,y_2) & =\frac12\left(x_1^2+y_1^2-x_2^2-y_2^2\right).
\end{align*}

\begin{proposition}\label{11}
The pair of momentum maps
\begin{equation}\label{pois}
\RR\stackrel{R}{\longleftarrow}(\CC^2,\om)\stackrel{J}{\longrightarrow} \su(2)^*=\RR^3
\end{equation}
is a dual pair.
\end{proposition}

\begin{proof}
The momentum maps for the commuting Hamiltonian actions of $\SU(2)$ and $S^1$ on $(\CC^2,\om)$ 
are equivariant, hence they form a pair of Poisson maps.
$R$ is obviously $\SU(2)$ invariant, so the dual pair property $(\ker TR)^\om=\ker TJ$  
follows from Proposition \ref{propo} if we show that $\SU(2)$ acts transitively  
on fibers of $R$.

To each element $\a\in\CC^2$ we associate a complex matrix 
$h_\a=\left[\begin{array}{ccc} a_1 & -\bar a_2 \\ a_2 & \bar a_1 \end{array}\right]$,
so the action of an element $g\in\SU(2)$ on $\CC^2$, $g\cdot\a=\b$, 
can be rewritten as matrix multiplication $g\cdot h_\a=h_\b$.
For any $r>0$, the fiber $R^{-1}(\frac {r^2}{2})$
is the 3-sphere $S_r^3$ of radius $r$.
We notice that $\frac{1}{r}h_\a\in\SU(2)$ for all $\a\in S_r^3$.
Given two elements in the same fiber, $\a,\b\in S_r^3$,
the matrix 
$$
g=h_\b h_\a^{-1}=\left(\frac1rh_\b\right)\left(\frac1r h_{\a}\right)^{-1}\in\SU(2)
$$ 
satisfies $g\cdot \a=\b$, hence $\SU(2)$ acts transitively on fibers of $R$.
\end{proof}

Because $X^2+Y^2+Z^2=R^2$, the momentum map $J=(X,Y,Z)$ maps the fibers of $R$, which are 3-spheres, 
into 2-spheres, coadjoint orbits of $\SU(2)$.
The restriction of $J$ to these 3-spheres is the {Hopf fibration}.
For this dual pair the symplectic leaf correspondence becomes
$\{c^2\}\mapsto J(R^{-1}(c^2))=S_{c}^2$.


\section{Poisson brackets on $\RR^3$}\label{sec4}

Vector fields $\v=( v_1, v_2, v_3)$ on $\RR^3$ with $ v_1, v_2, v_3\in\F(\RR^3)$ are in 1-1 correspondence with
bivector fields on $\RR^3$
\[
\pi_\v= v_1\pa_y\wedge\pa_z+ v_2\pa_z\wedge\pa_x+ v_3\pa_x\wedge\pa_y.
\]
The following are necessary and sufficient conditions for the bivector field $\pi_\v$ to be Poisson:
\begin{enumerate}
\item $ v_1\left(\pa_y v_3-\pa_z v_2\right)
+ v_2\left(\pa_z v_1-\pa_x v_3\right)
+ v_3\left(\pa_x v_2-\pa_y v_1\right)=0$.
\item $\v^\flat\wedge d(\v^\flat)=0$, where $\v^\flat= v_1 dx+ v_2 dy+  v_3 dz$.
\item The distribution $\v^\perp$ on $\RR^3$ is integrable.
\end{enumerate}
Under these circumstances the Hamiltonian vector field with Hamiltonian function $H$ 
on the Poisson manifold $(\RR^3,\pi_\v)$ is $X_H=\v\x\nabla H$,
with $\x$ denoting the usual vector product on $\RR^3$,
so the Poisson bracket on $\RR^3$ associated to the bivector field $\pi_\v$  can be written as
\[
\{F,G\}_\v=\v\cdot(\nabla F\x\nabla G).
\]
All Hamiltonian vector fields are orthogonal to $\v$, hence the symplectic leaves of the Poisson structure $\pi_\v$
are leaves of the integrable distribution $\v^\perp$.

\bigskip

The equivalent conditions 1., 2., and 3. are satisfied for gradient vector fields $\v=\nabla C$ with $C\in\F(\RR^3)$. 
The associated Poisson bracket is the {\bf Nambu bracket} 
\[
\{F,G\}_{\nabla C}=\nabla C\cdot(\nabla F\x\nabla G)=Jac(C,F,G),
\]
where $Jac$ denotes the Jacobian determinant.
The function $C$ is a Casimir  and the symplectic leaves are the level surfaces $C=\text{constant}$.
A similar result holds in a more general setting:

\begin{proposition}\label{rem}
The vector field $\v=f\nabla C$, where $f$ is a nonvanishing function on $\RR^3$,
determines a Poisson structure $\pi_\v$ on $\RR^3$ with symplectic leaves the level surfaces of the function $C$.
\end{proposition}

\begin{proof}
From the three equivalent conditions, the third one is the easiest to check:
the distribution $\v^\perp$ coincides with the orthogonal distribution to the gradient vector field of $C$,
hence it is integrable.
\end{proof}


\section{Kummer shapes as symplectic leaves}\label{sec5}

\begin{figure}
\begin{center}
\includegraphics[width=29mm]{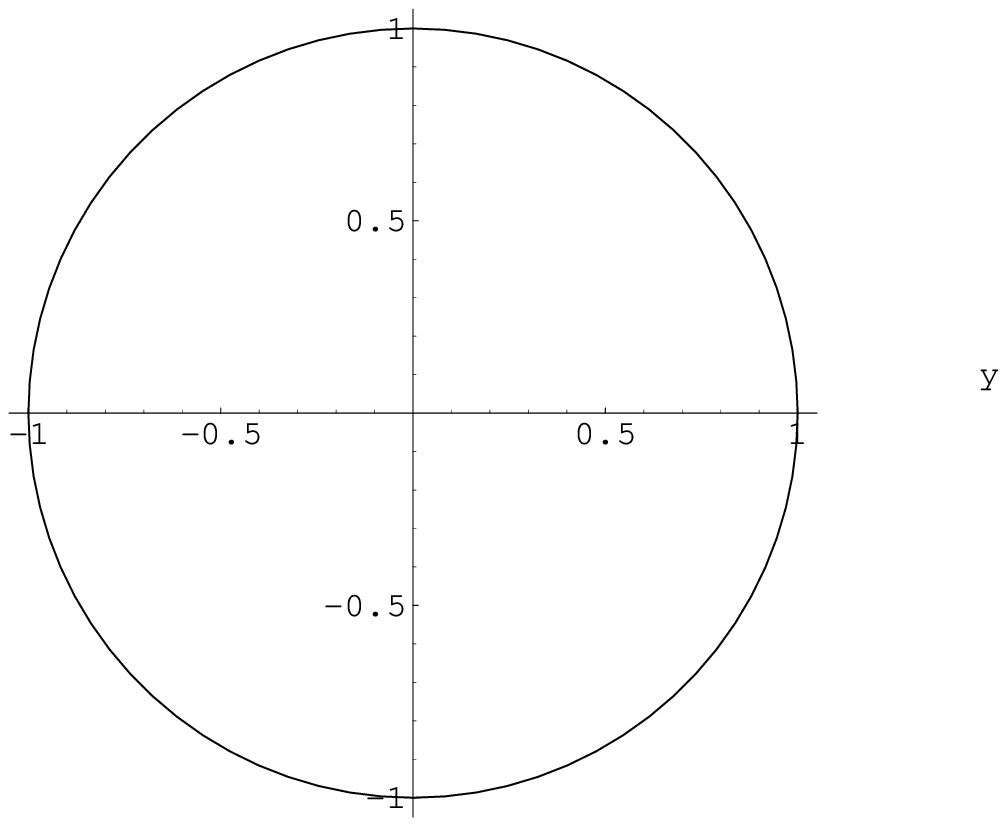}
\includegraphics[width=29mm]{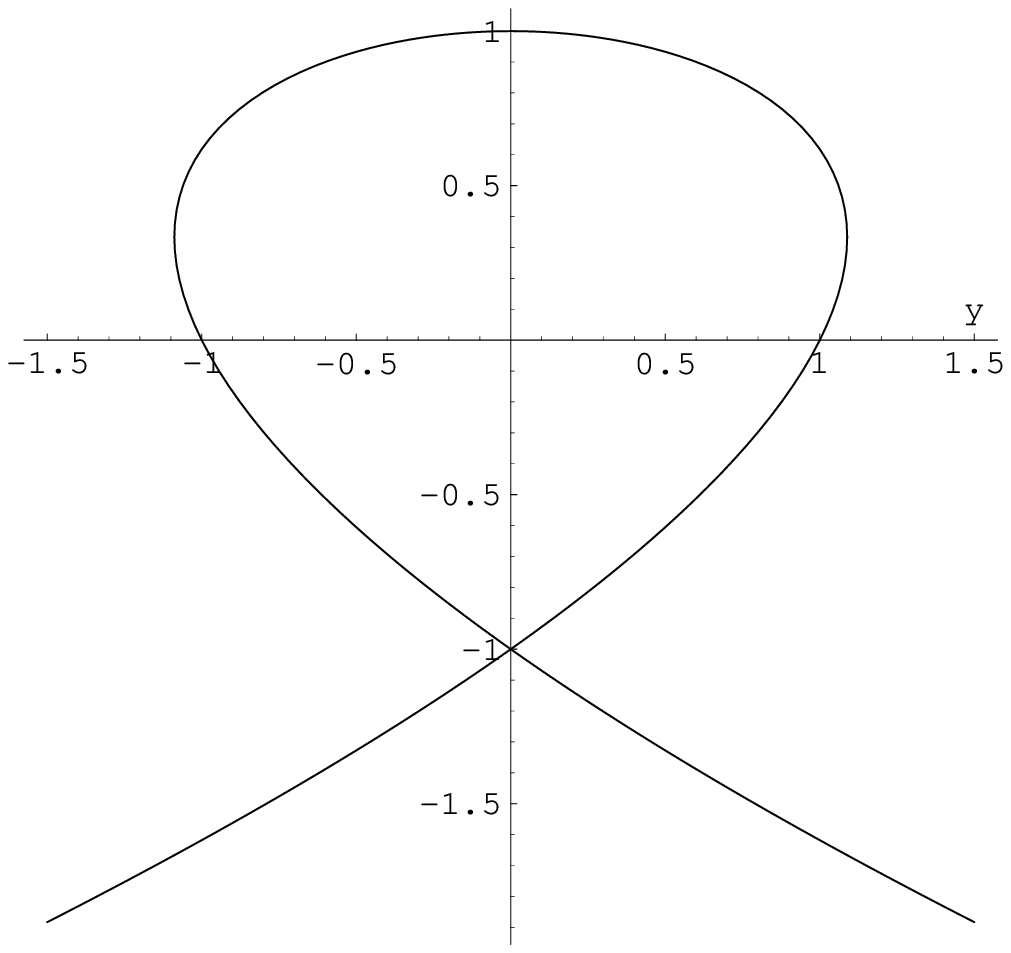}
\includegraphics[width=31mm]{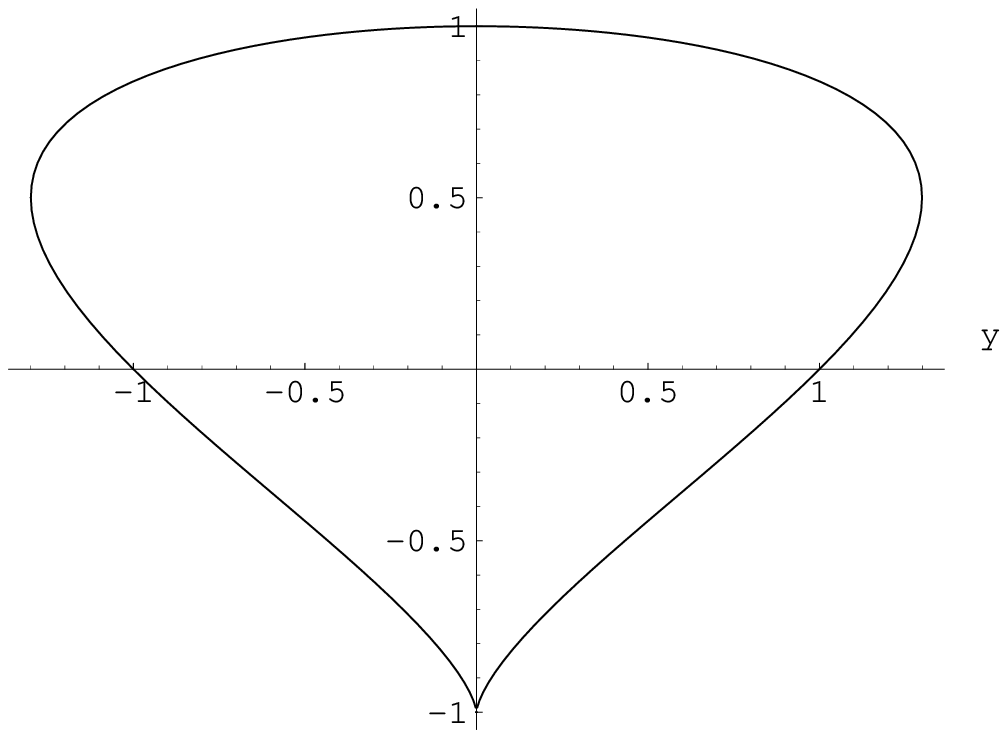}
\includegraphics[width=30mm]{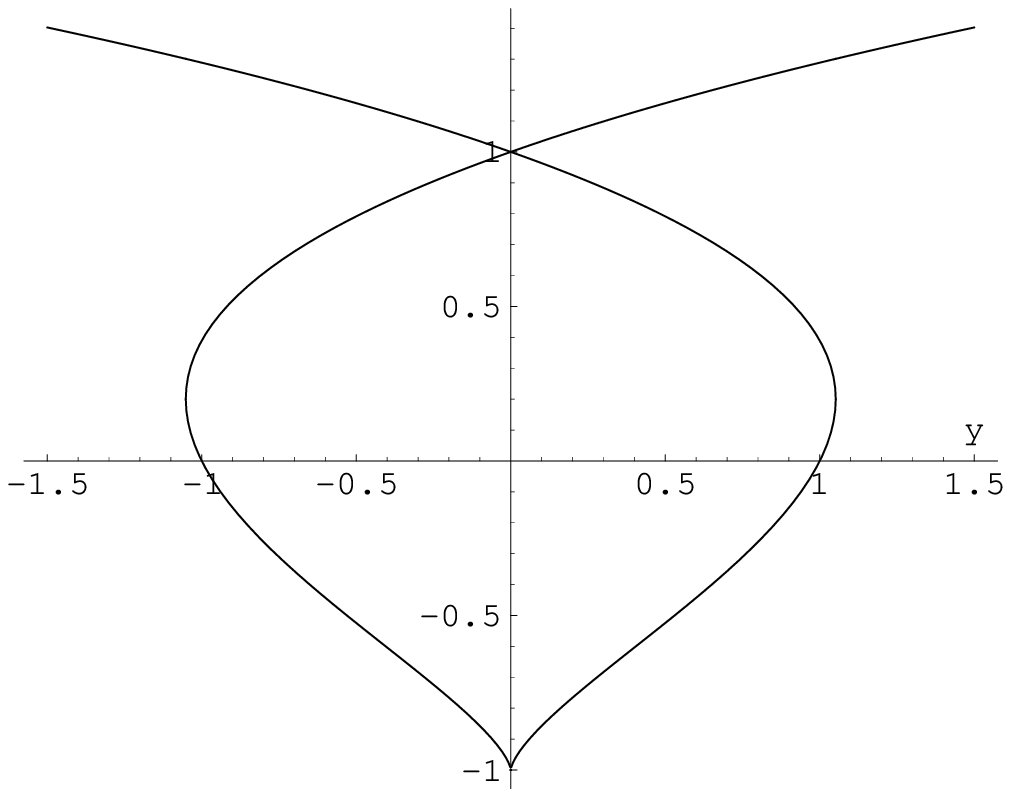}
\includegraphics[width=31mm]{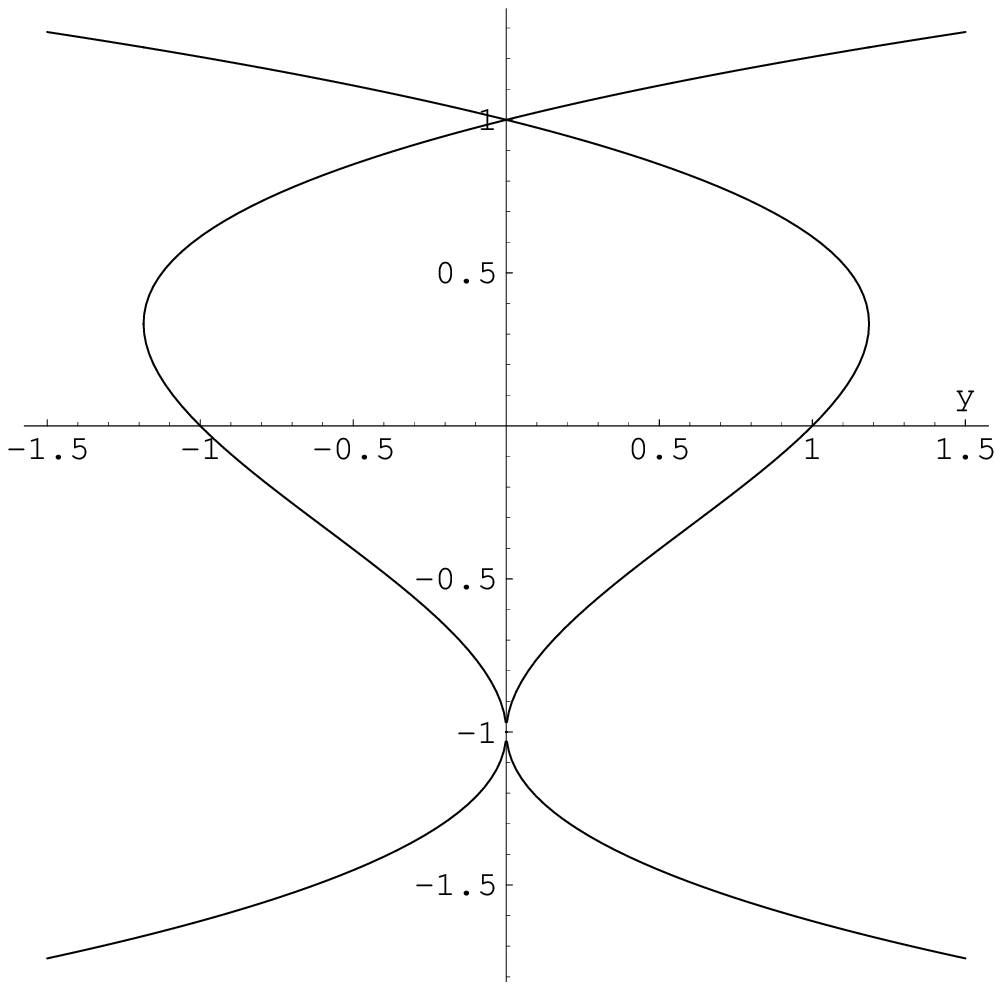}
\end{center}
\caption{The curve $y^2=\left(\frac{1+z}{n}\right)^m\left(\frac{1-z}{m}\right)^n$ for $(m,n)$ equal to (1,1), (2,1), (3,1), (3,2) and (4,2)}
\label{fig8}
\end{figure}

\begin{figure}
\begin{center}
\includegraphics[width=26mm]{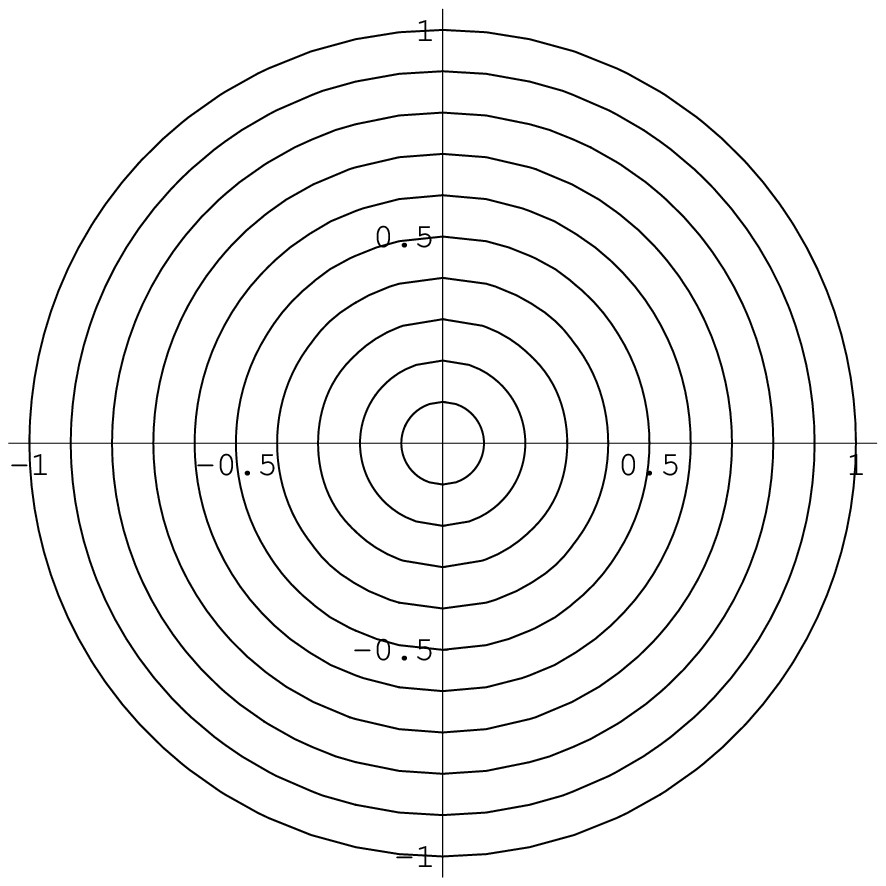}
\includegraphics[width=29mm]{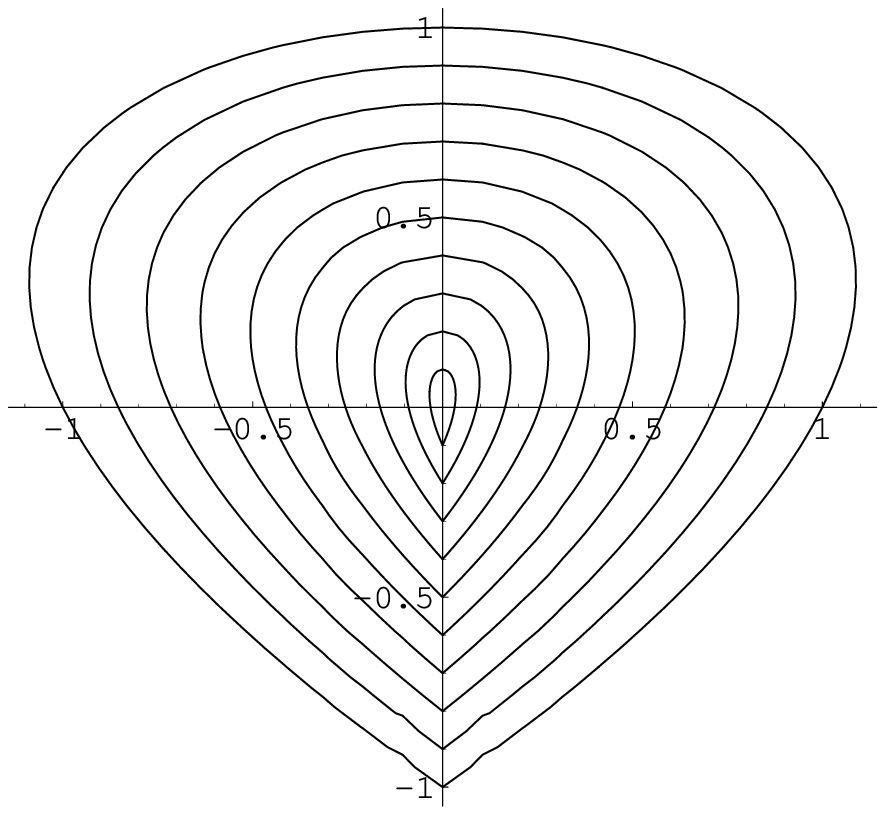}
\includegraphics[width=32mm]{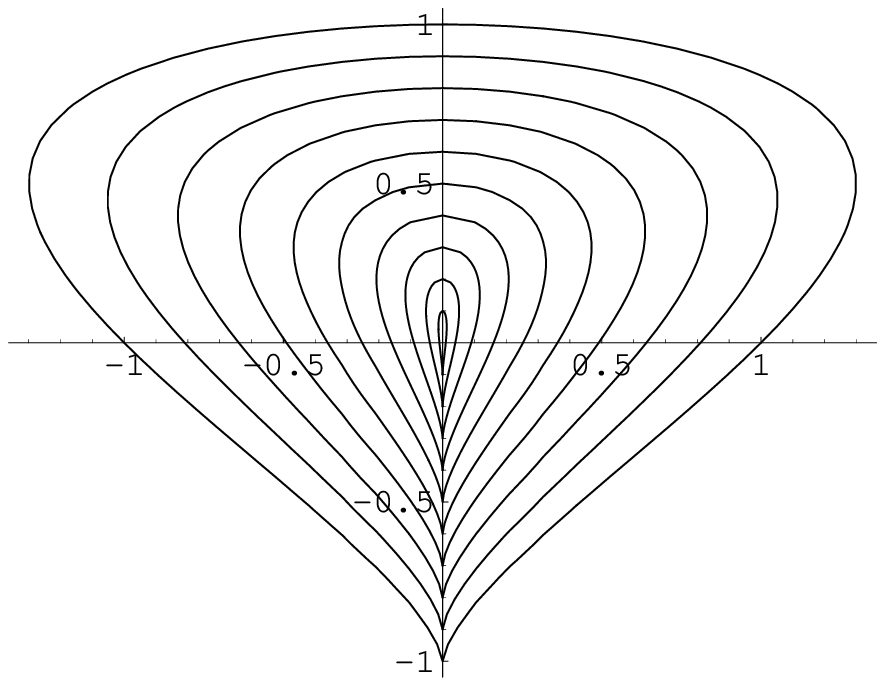}
\includegraphics[width=28mm]{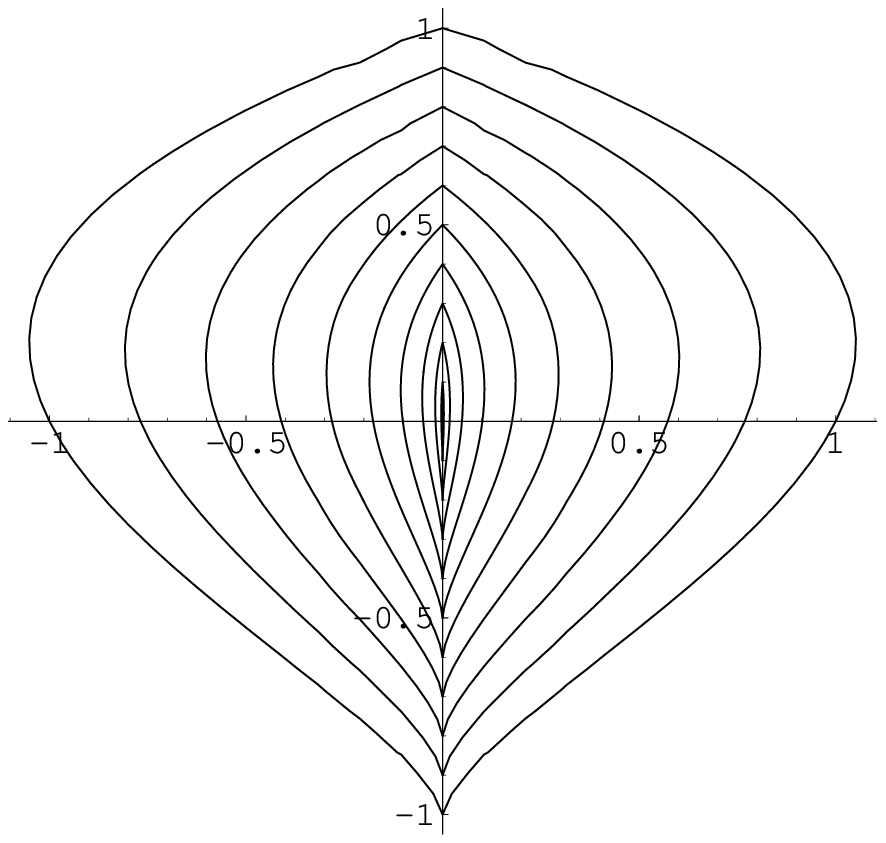}
\includegraphics[width=31mm]{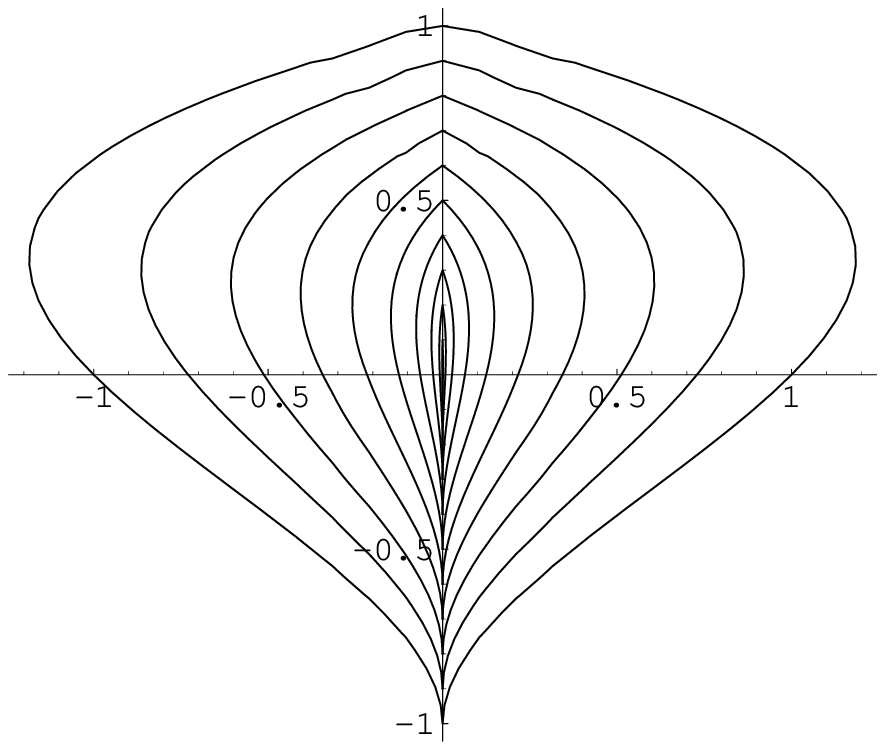}
\end{center}
\caption{Curves generating Kummer shapes for 1:1, 2:1, 3:1, 3:2 and 4:2 resonance}
\label{fig1}
\end{figure}

The {\bf Kummer shapes} in $n:m$ resonance, $n,m>0$, are the bounded surfaces defined by 
the equation \cite{Ku1986}
\begin{equation}\label{k}
x^2+y^2-\left(\frac{c+z}{n}\right)^m\left(\frac{c-z}{m}\right)^n=0,
\quad |z|<c,
\end{equation}
where $c$ is a positive constant (see Figure \ref{fig1}).
They are obtained by rotating around the $z$ axis the algebraic curve 
(see Figure \ref{fig8})
$$
y^2=\left(\frac{c+z}{n}\right)^m\left(\frac{c-z}{m}\right)^n,\quad |z|<c.
$$
Let $\Phi\in\F(\RR^4)$ be given by
\begin{equation}\label{fi}
\Phi(x,y,z,r)=x^2+y^2-\left(\frac{r+z}{n}\right)^m\left(\frac{r-z}{m}\right)^n.
\end{equation}
One can obtain the Kummer shapes also by slicing with hyperplanes $r=c$ of $\RR^4$ 
that part of the hypersurface $\Phi=0$ included in the intersection of the halfspaces $z<r$ and $z>-r$.

\begin{lemma}\label{hard}
The Kummer shapes can be expressed as level sets of a smooth function $C$ defined 
on $\RR^3$ with the $z$ axis removed. 
\end{lemma}

\begin{proof}
We must show there exists a smooth function $C:\RR^3\setminus Oz\to\RR$ which satisfies
\begin{equation}\label{imp}
\Phi(x,y,z,C(x,y,z))=0,\quad |z|<C(x,y,z).
\end{equation}
In other words the hypersurface $\Phi=0$ of $\RR^4$ coincides with the graph of the function $C$,
on the intersection of the halfspaces $z<r$ and $z>-r$.

To prove this, we use rotational symmetry in $(x,y)$ about the $z$-axis to reduce the problem to 
proving the existence and uniqueness of a smooth function $c:(0,\infty)\x\RR\to\RR$
such that 
\[
y^2=\left(\frac{c(y,z)+z}{n}\right)^m\left(\frac{c(y,z)-z}{m}\right)^n,\quad |z|<c(y,z).
\]
Then $C(x,y,z)=c(\sqrt{x^2+y^2},z)$ is a smooth function defined on $\RR^3\setminus Oz$ whose level sets are the Kummer shapes.

The polynomial function $p(r)=\left(\frac{r+z}{n}\right)^m\left(\frac{r-z}{m}\right)^n-y^2$,
with coefficients smoothly depending on $(y,z)\in(0,\infty)\x\RR$,
has at least one zero in the interval $(|z|,\infty)$
because $p(|z|)=-y^2<0$ and $\lim_{r\to\infty}p(r)=+\infty$.
But $p$ is a monotone increasing function on  $(|z|,\infty)$, so
there is a unique zero of $p$ in the interval $(|z|,\infty)$, denoted by $c(y,z)$.
This gives the smooth function $c$ on $(0,\infty)\x\RR$ we were seeking.
\end{proof}

From the implicit identity \eqref{imp} we deduce that the gradient vector field $\nabla C$ can be written as
\begin{equation*}
\nabla C=\left(-\frac{1}{\pa_r\Phi}{\nabla_{(x,y,z)}\Phi}\right)\Big|_{r=C}.
\end{equation*}
It follows that the vector field $\v$, defined on $\RR^3\setminus Oz$ by
\begin{align}\label{v}
\v&:=\nabla_{(x,y,z)}\Phi\big|_{r=C}=\left(2x,2y,-(x^2+y^2)\left(\frac{m}{C(x,y,z)+z}-\frac{n}{C(x,y,z)-z}\right)\right),
\end{align}
is of the form $\v=f\nabla C$, where $f$ is the  nonvanishing function 
\begin{equation}\label{f}
f=-{\pa_r\Phi}|_{r=C}.
\end{equation}

\begin{proposition}
The Kummer shapes \eqref{k}, with the singular points $(0,0,\pm c)$ removed, 
are symplectic leaves of the Poisson manifold $(\RR^3\setminus Oz,\pi_\v)$ 
associated to the vector field $\v$ given by \eqref{v}.
\end{proposition}

\begin{proof}
We know from \eqref{v} that $\v=f\nabla C$,
so by proposition \ref{rem} the bivector field $\pi_\v$ is a Poisson bivector field. Its symplectic leaves 
are the surfaces $C=\text{constant}$, \ie the Kummer shapes.
\end{proof}



\section{Poisson maps for $n:m$ resonance}\label{sec6}

Let $n$ and $m$ be non-zero natural numbers.
The action 
\begin{equation}\label{act}
z\cdot(a_1,a_2)=(z^na_1,z^ma_2),\quad z\in S^1\subset\CC
\end{equation}
of the circle $S^1$ on $\CC^2$, with the opposite $\om$ of the canonical symplectic form,
\begin{equation*}
\omega=-\frac i2 (da_1\wedge d\bar a_1+da_2\wedge d\bar a_2)=-dx_1\wedge dy_1-dx_2\wedge dy_2
\end{equation*}
is Hamiltonian
with infinitesimal action $(a_1,a_2)\mapsto (ina_1,ima_2)$.
The associated momentum map 
\begin{equation}\label{r+}
R:\CC^2\to\RR,\quad R(\a)=\frac n2|a_1|^2+\frac m2|a_2|^2
\end{equation}
is equivariant, which implies that $R$ is a Poisson map.

\bigskip

Let $X,Y,Z$ be the functions on $\CC^2$ uniquely defined by the identities
\begin{equation}\label{xyz}
X(\a)-iY(\a)=a_1^m\bar a_2^{n}\text{ and } Z=\frac n2|a_1|^2-\frac m2|a_2|^2.
\end{equation}
An easy computation reveals that 
\[
X^2+Y^2
=\left(\frac{R+Z}{n}\right)^m\left(\frac{R-Z}{m}\right)^n.
\]
This can be written as $\Phi\o(X,Y,Z,R)=0$ on $\CC^2$, with $\Phi$ the function \eqref{fi}, which means that $C\o(X,Y,Z)=R$ on $(\CC\setminus\{0\})^2$. Here we have to restrict the functions $X,Y,Z,R$ to 
$(\CC\setminus\{0\})^2$ because $C$ is not defined on the $z$ axis. 

\begin{proposition}\label{p2}
The map $\Pi=(X,Y,Z):(\CC\setminus\{0\})^2\to\RR^3\setminus Oz$ is a Poisson map with respect to $\om$, the 
opposite of the canonical symplectic form on $(\CC\setminus\{0\})^2$ and the Poisson bivector field $\pi_{mn\v}$ on $\RR^3\setminus Oz$,
with vector field $\v$ defined by \eqref{v}.
\end{proposition}

\begin{proof}
The following Poisson brackets on the symplectic manifold $(\CC^2,\om)$ are computed in \cite{Holm}:
\begin{align*}
\{Y,Z\}&=2mnX\\
\{Z,X\}&=2mnY\\
\{X,Y\}&=-\,mn(X^2+Y^2)\left(\frac{m}{R+Z}-\frac{n}{R-Z}\right).
\end{align*}
Knowing that
\begin{align*}
\pi_{mn\v}=2mnx\pa_y\wedge\pa_z+2mny\pa_z\wedge\pa_x
-mn(x^2+y^2)\left(\frac{m}{C(x,y,z)+z}-\frac{n}{C(x,y,z)-z}\right)\pa_x\wedge\pa_y,
\end{align*}
the result follows from the functional identity $C\o(X,Y,Z)=R$.
\end{proof}

One may also verify that $\Pi$ is a surjective submersion.

\begin{remark}
{\rm For $n=m=1$ there are no singularities along the $z$ axis, so one obtains the Poisson structure 
$$
2x\pa_y\wedge\pa_z+2y\pa_z\wedge\pa_x+2z\pa_x\wedge\pa_y
$$
on all of $\RR^3$. This is isomorphic to the Lie-Poisson structure on $\su(2)^*$, the dual of the Lie algebra of $\SU(2)$.
The Kummer shapes are spheres: the coadjoint orbits of $\SU(2)$.
Moreover, in this case the map $\Pi$ becomes the equivariant momentum map $J:\CC^2\to\su(2)^*$ from \eqref{j_can},
for the canonical Hamiltonian $\SU(2)$-action on $\CC^2$.}
\end{remark}


\section{The $n:m$ resonance dual pair}\label{sec7}

We saw in Proposition \ref{11} that
the dual pair for 1:1 resonance is the pair $(R,J)$ 
of  momentum maps associated to the natural commuting Hamiltonian actions 
of $S^1$ and ${\SU(2)}$ on  $\CC^2$ with the opposite $\om$ of the canonical symplectic form:
\[
{\RR\stackrel{R}{\longleftarrow}(\CC^2,\om)
\stackrel{J}{\longrightarrow} \su(2)^*=\RR^3}.
\] 
There is a dual pair also for general $n:m$ resonance,
but it is a dual pair of Poisson maps, rather than momentum maps.

On $\RR^3\setminus Oz$ we consider the Poisson bivector field $\pi_{mn\v}$,
with $\v$ the vector field \eqref{v}.

\begin{theorem}\label{t1}
The pair of Poisson maps
\[
{\RR\stackrel{R}{\longleftarrow}((\CC\setminus\{0\})^2,\om)
\stackrel{\Pi}{\longrightarrow}(\RR^3\setminus Oz,\pi_{mn\v})}
\] 
is a dual pair for all pairs $(m,n)$ of nonzero natural numbers.
\end{theorem}

\begin{proof}
We know already from the previous section that both $R$ and $\Pi$ are Poisson maps,
we only have to check the dual pair property 
\begin{equation}\label{dpbis}
\ker T_\a R=(\ker T_\a\Pi)^\om,\quad\forall\a\in (\CC\setminus\{0\})^2.
\end{equation}

The symplectic form $\om$ and the canonical Riemannian metric $g$ on $\CC^2$ introduced in Section 3 are related by
$\om(\a,\b)=g(\a,i\b)$, so the symplectic and Riemannian orthogonals to a real vector subspace $V\subset\CC^2$ 
are also related: $V^\perp=(iV)^\om$. If the vector subspace $V$ is generated by the vector $\a=(a_1,a_2)\in\CC^2$,
then $V^\om=\a^\om$ and $V^\perp=\a^\perp$. Thus we get 
\begin{equation}\label{one}
\ker T_{\a}R=(na_1,ma_2)^\perp=(nia_1,mia_2)^\om.
\end{equation}

Using the expression \eqref{xyz} of the functions $X,Y,Z$, it is not hard to verify that 
\[
(nia_1,mia_2)\in\ker T_{\a}\Pi=\ker T_{\a}X\cap\ker T_{\a}Y\cap\ker T_{\a}Z
\]
We check it here for the function $X(\a)=\frac12(a_1^m\bar a_2^{n}+\bar a_1^{m}a_2^n)$,
the computations being similar for $Y$ and $Z$ from \eqref{xyz}:
\begin{align*}
T_{\a}X.(nia_1,mia_2)
&=\frac12\big(ma_1^{m-1}\bar a_2^{n}(nia_1)
+na_1^m\bar a_2^{n-1}(-mi\bar a_2)+m\bar a_1^{m-1}a_2^n(-ni\bar a_1)
+n\bar a_1^{m}a_2^{n-1}(mia_2)\big)=0
\end{align*}
implies $(nia_1,mia_2)\in\ker T_{\a}X$.
The kernel of $T_{\a}\Pi$ is 1-dimensional ($\Pi$ is a submersion), 
so it must be generated by the nonzero vector $(nia_1,mia_2)$.
We get
\[
(\ker T_{\a}\Pi)^\om=(nia_1,mia_2)^\om,
\]
which, together with \eqref{one}, ensures the dual pair property \eqref{dpbis}.
\end{proof}

The symplectic leaf correspondence theorem for dual pairs, 
applied to the $n:m$ resonance, says that, for each $c>0$, 
the symplectic leaf $\{c\}$ of $\RR$ corresponds to the symplectic leaf
$\Pi(R^{-1}(c))$ of $\RR^3$, \ie to the Kummer surface 
$C(x,y,z)=c$, because $C\o\Pi=R$.


\section{The $1:-1$ resonance dual pair}\label{sec8}

In this section we give an alternative approach to \cite{Iwai} for the $1:-1$ resonance dual pair.
The Lie group $\U(1,1)$ of complex $2\x 2$ matrices preserving the Hermitian inner product 
\begin{equation}\label{herminus}
\langle\a,\b\rangle_-=a_1\bar b_1-a_2\bar b_2 \ \text{ on }\CC^2.
\end{equation}
has a 1-dimensional center, 
isomorphic to $S^1$, and a normal subgroup $\SU(1,1)$ consisting of complex matrices with determinant 1:
$$
\SU(1,1)=\left\{\left[\begin{array}{ccc} \al & \be \\ \bar\be & \bar\al \end{array}\right]\text{ with }|\al|^2-|\be|^2=1\right\}.
$$
We endow $\CC^2$ with the symplectic form 
\begin{equation}\label{omega}
\omega_-=-\frac i2(da_1\wedge d\bar a_1-da_2\wedge d\bar a_2)=-dx_1\wedge dy_1+dx_2\wedge dy_2,
\end{equation}
the imaginary part of the inner product \eqref{herminus}.
The natural action of the group $\U(1,1)$ by multiplication on $\CC^2$ is Hamiltonian,
with $\U(1,1)$-equivariant momentum map 
\begin{equation}\label{canon-}
\bar J_-:\CC^2\to\u(1,1)^*,\quad \langle\bar J_-(\a),\xi\rangle_{\u(1,1)}=\frac12(-\om_-)(\xi(\a),\a)
=\frac i2 \langle\a,\xi(\a)\rangle_-.
\end{equation}

The map 
\begin{equation}\label{r-}
R_-:\CC^2\to\RR,\quad R_-=\frac 12|a_1|^2-\frac 12|a_2|^2.
\end{equation}
is the momentum map for the natural $S^1$-action $z\cdot(a_1,a_2)=(za_1,za_2)$
on $(\CC^2,\om_-)$. This is the action of the center of $\U(1,1)$.
There is a linear isomorphism between $\RR^3$ and the Lie algebra $\su(1,1)$ given by
$$
u\in\RR^3\mapsto 
\left[\begin{array}{ccc} iu_3 & iu_1+u_2 \\- iu_1+u_2 & -iu_3 \end{array}\right]\in\su(1,1).
$$  
With this identification, from the expression \eqref{canon-} of $\bar J_-$ we deduce the following expression
of the momentum map for the $\SU(1,1)$ action: 
\begin{equation}\label{j-}
J_-:\CC^2\to \su(1,1)^*=\RR^3,\quad J_-(\a)=\left({\rm Re}(a_1 \bar a_2) , -\,{\rm Im}(a_1 \bar a_2) ,
-\left(\frac 12|a_1|^2+\frac 12|a_2|^2\right) \right).
\end{equation}
Denoting by $(X,Y,Z_-)$ the three components of the momentum map $J_-$, 
we get that $X^2+Y^2-Z_-^2=R_-^2$. 

\begin{proposition}
The pair of momentum maps \eqref{r-} and \eqref{j-} for the commuting actions
of $S^1$ and $\SU(1,1)$ on $(\CC^2,\om_-)$ 
\[
{\RR\stackrel{R_-}{\longleftarrow}(\CC^2,\om_-)
\stackrel{J_-}{\longrightarrow} \su(1,1)^*=\RR^3}
\]  
is a dual pair.
\end{proposition}

The proof is similar to that of Proposition \ref{11}. It uses the fact that
the action of an element $g\in\SU(1,1)$ on $\CC^2$, $g\cdot\a=\b$, 
can be rewritten as matrix multiplication $g\cdot k_\a=k_\b$,
where 
$k_\a=\left[\begin{array}{ccc} a_1 & \bar a_2 \\ a_2 & \bar a_1 \end{array}\right]$.
Given two elements $\a,\b$ in the same fiber of $R_-$, 
the matrix $g=k_\b k_\a^{-1}\in\SU(1,1)$ 
satisfies $g\cdot \a=\b$, hence $\SU(1,1)$ acts transitively on fibers of $R_-$.

\bigskip

The momentum map $J_-$ maps the fibers of $R_-$, which are 3-hyperboloids, 
into 2-hyperboloids, coadjoint orbits of $\SU(1,1)$.
The restriction of $J_-$ to these 3-hyperboloids is the hyperbolic Hopf fibration.


\section{The $n:-m$ resonance dual pair}\label{sec9}

\begin{figure}
\begin{center}
\includegraphics[width=29mm]{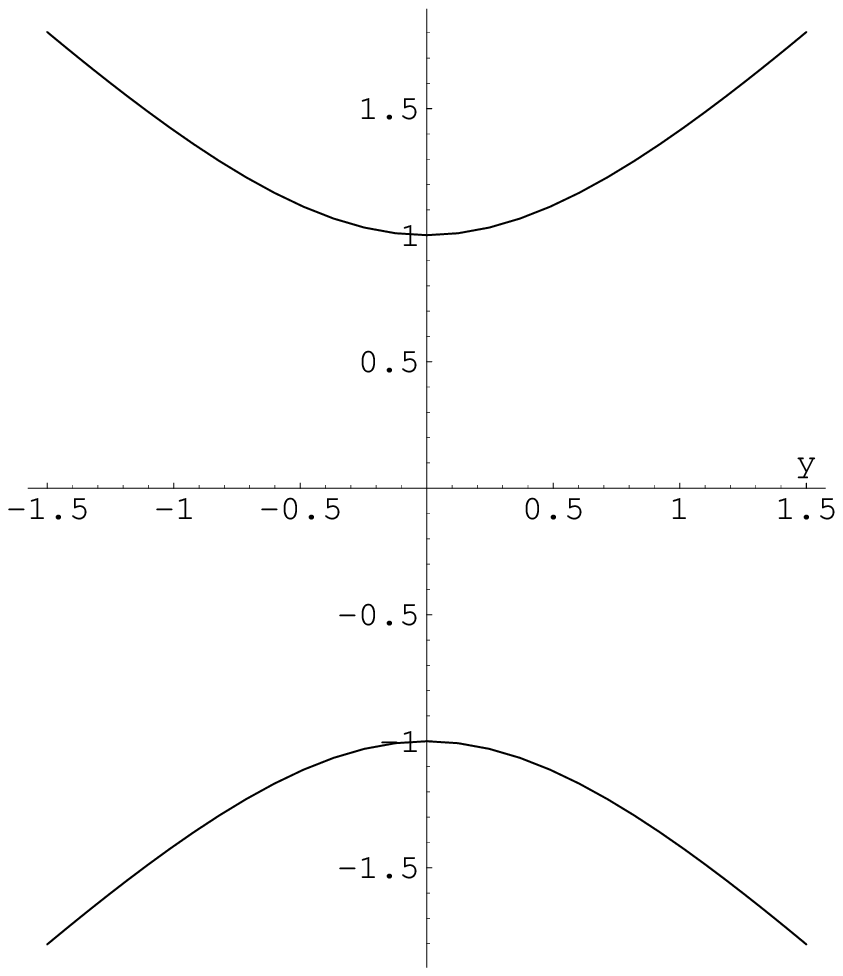}
\includegraphics[width=29mm]{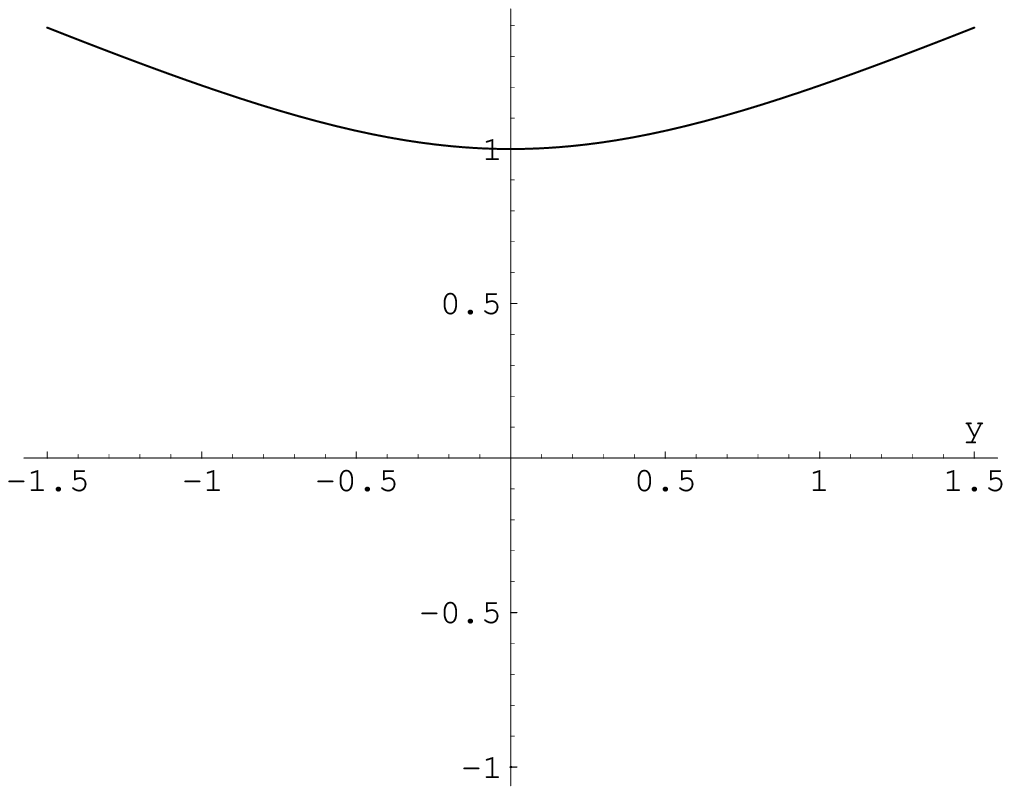}
\includegraphics[width=31mm]{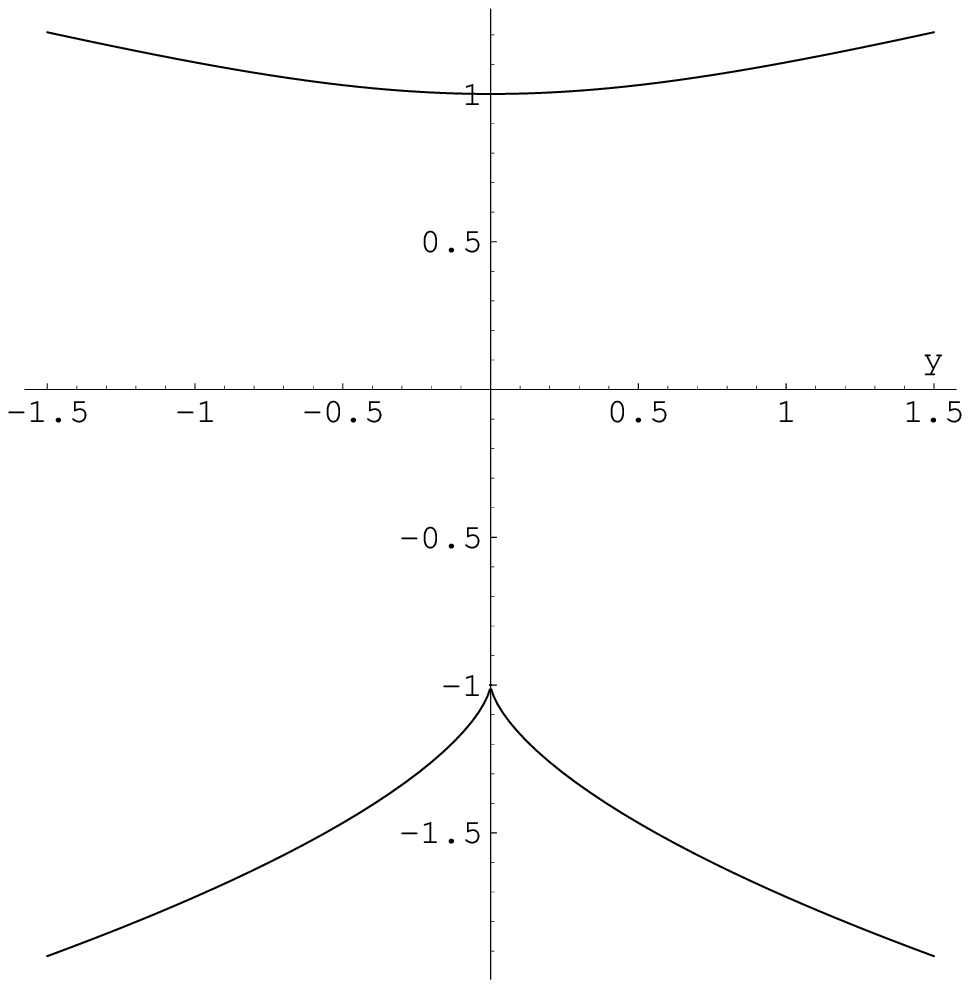}
\includegraphics[width=30mm]{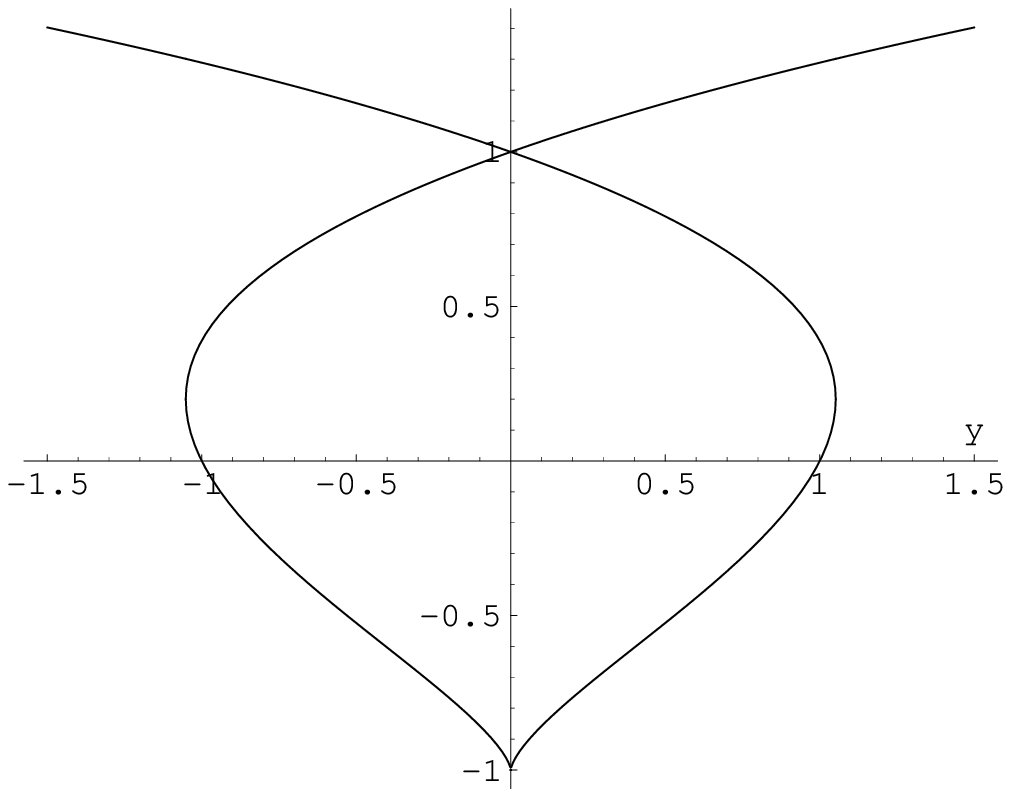}
\includegraphics[width=31mm]{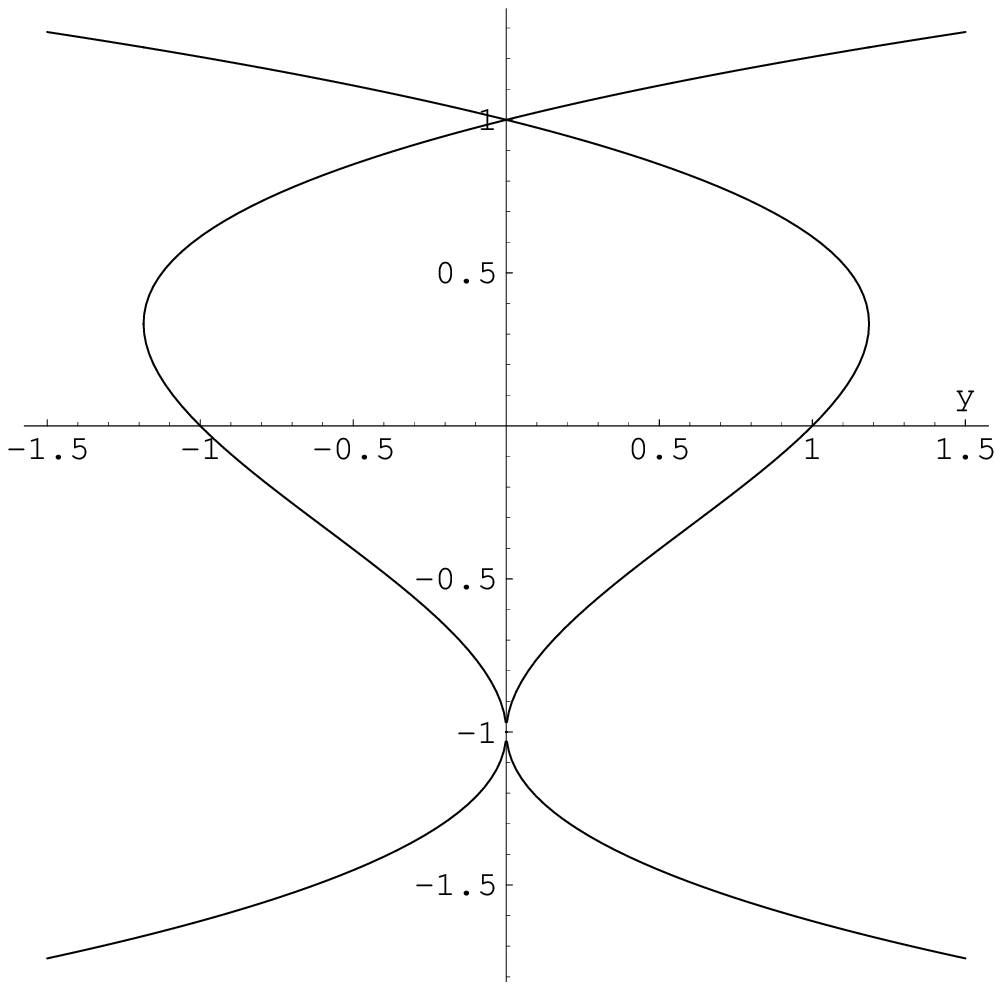}
\end{center}
\caption{The curve $y^2=\left(\frac{z+1}{n}\right)^m\left(\frac{z-1}{m}\right)^n$ for $(m,n)$ equal to (1,1), (2,1), (3,1), (3,2) and (4,2)}
\label{fig7}
\end{figure}

In this section we build a dual pair for the more general $n:-m$ resonance,
a dual pair in which the second map is not a momentum map, but only a Poisson map.
The first map is 
\[
R_-:\CC^2\to\RR,\quad R_-(\a)=\frac n2|a_1|^2-\frac m2|a_2|^2,
\]
the equivariant momentum map for the $S^1$-action \eqref{act} on $(\CC^2,\om_-)$.
The second map is $\Pi_-=(X,Y,Z_-)$, where 
\begin{equation}\label{nega}
Z_-:\CC^2\to\RR,\quad Z_-(\a)=\frac n2|a_1|^2+\frac m2|a_2|^2,
\end{equation}
and, as for the $n:m$ resonance, the functions $X,Y:\CC^2\to\RR$  are defined by
the identity
\[
X(\a)-iY(\a)=a_1^m\bar a_2^{n}.
\]

\begin{figure}
\begin{center}
\includegraphics[width=26mm]{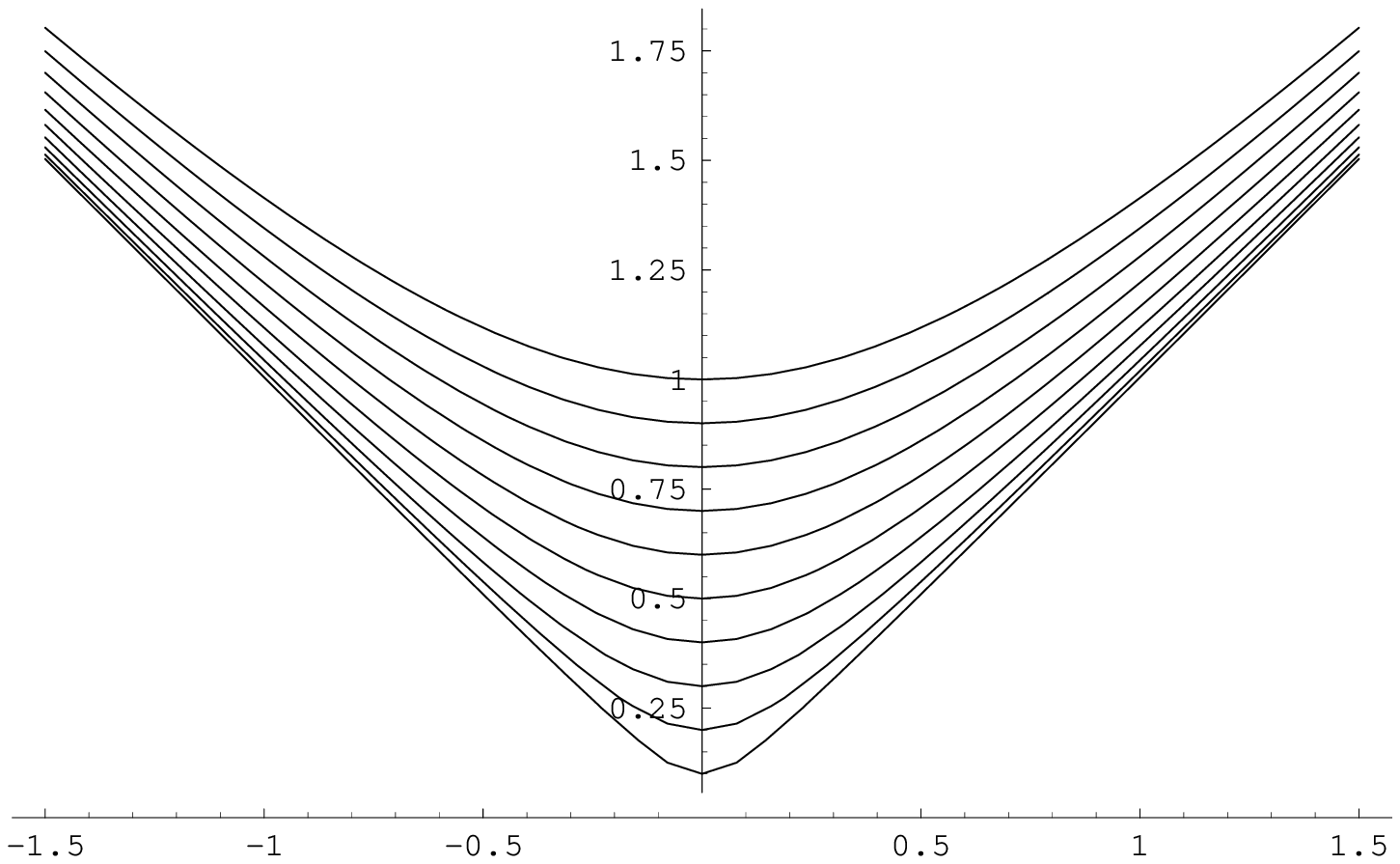}
\includegraphics[width=29mm]{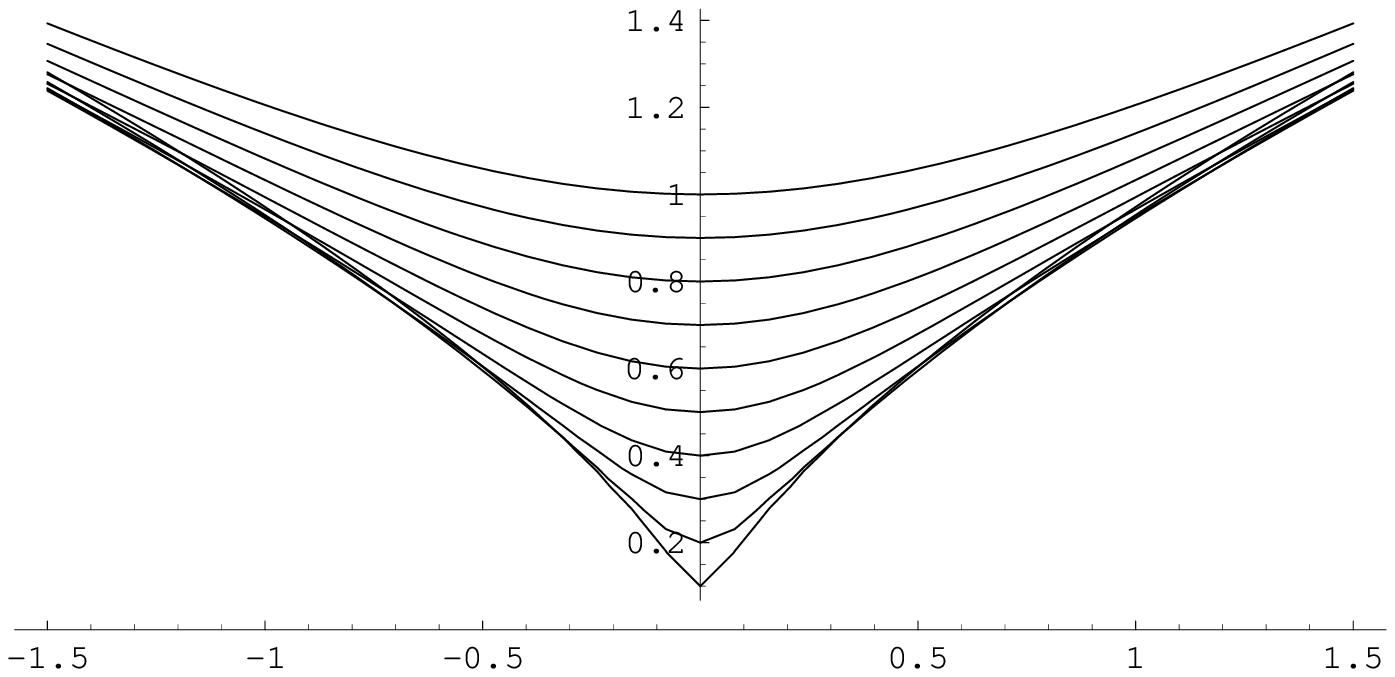}
\includegraphics[width=34mm]{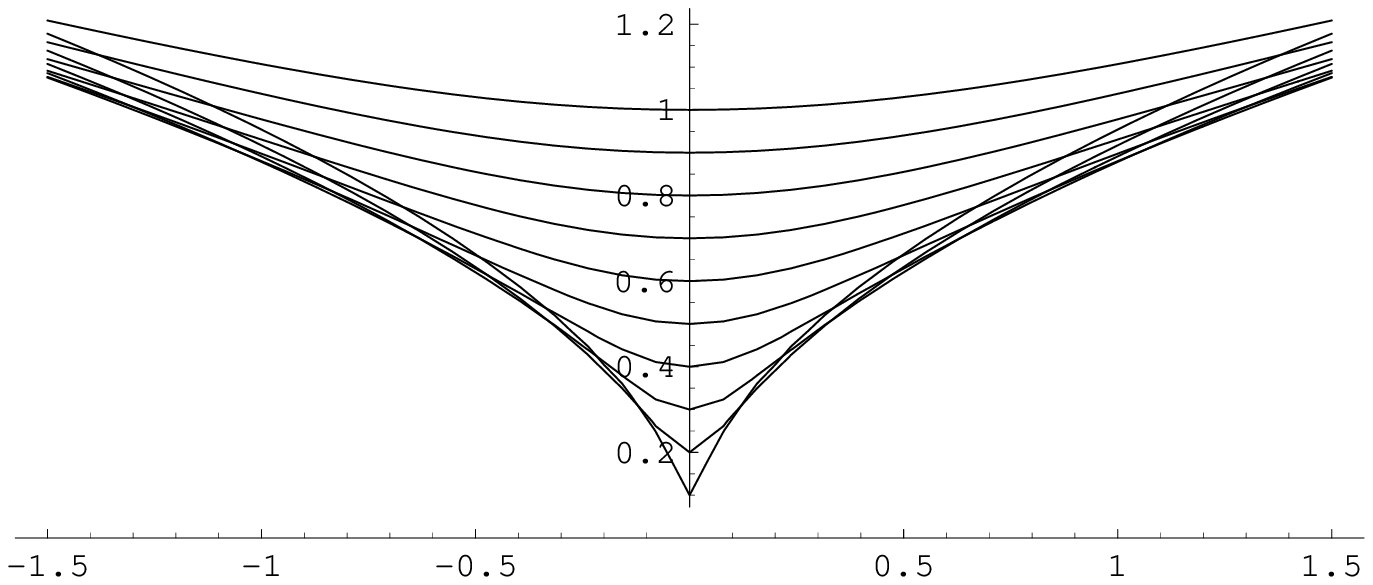}
\includegraphics[width=28mm]{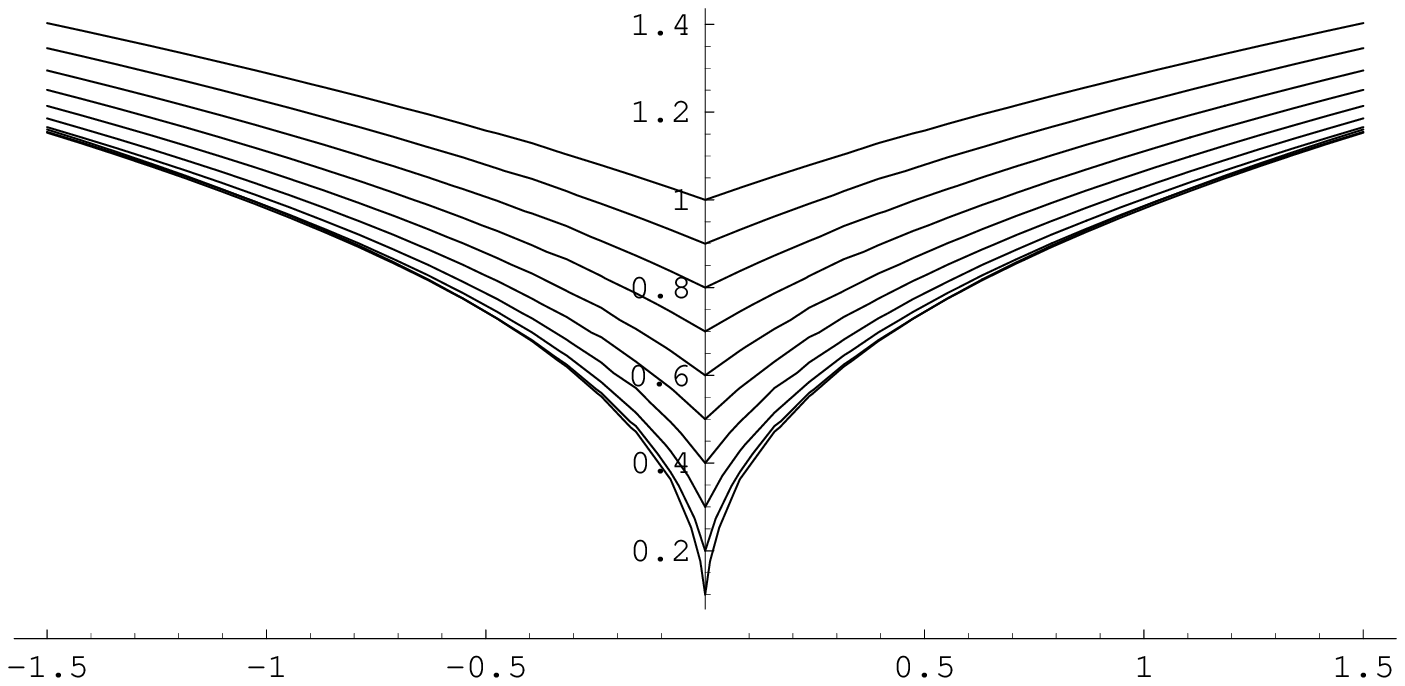}
\includegraphics[width=31mm]{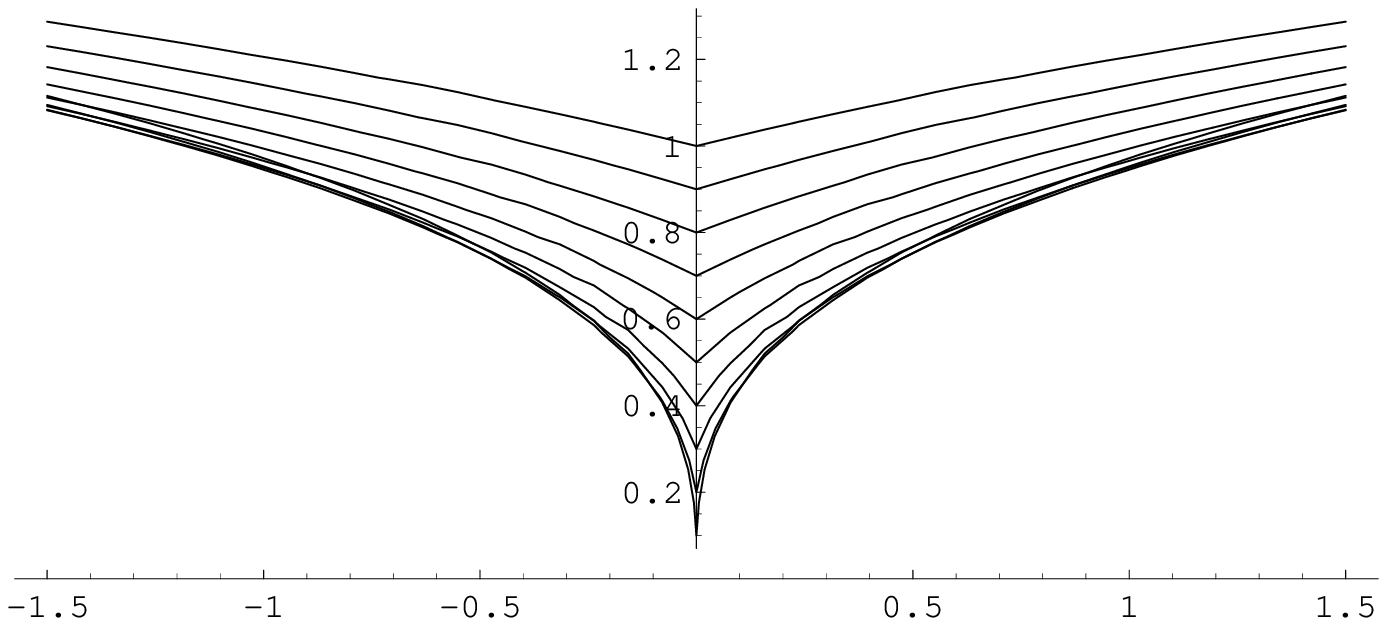}
\end{center}
\caption{Curves generating upper Kummer shapes for 1:-1, 2:-1, 3:-1, 3:-2 and 4:-2 resonance}
\label{fig3}
\end{figure}

\begin{figure}
\begin{center}
\includegraphics[width=27mm]{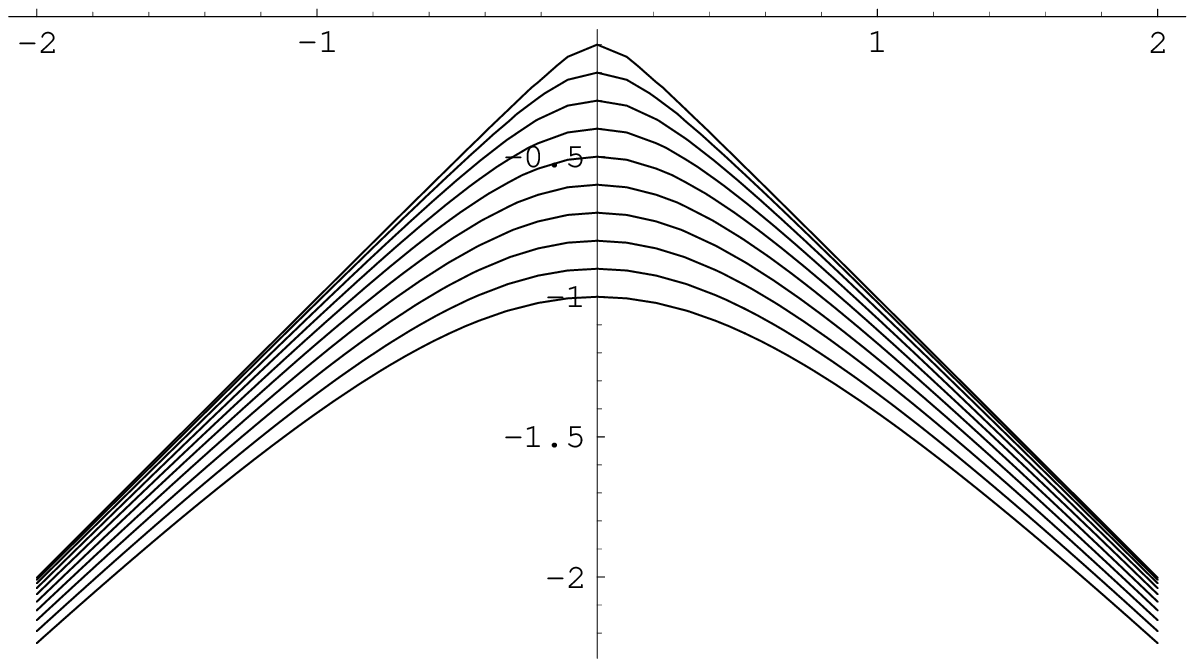}
\quad\quad\quad\quad\quad\quad\quad\quad\quad
\includegraphics[width=30mm]{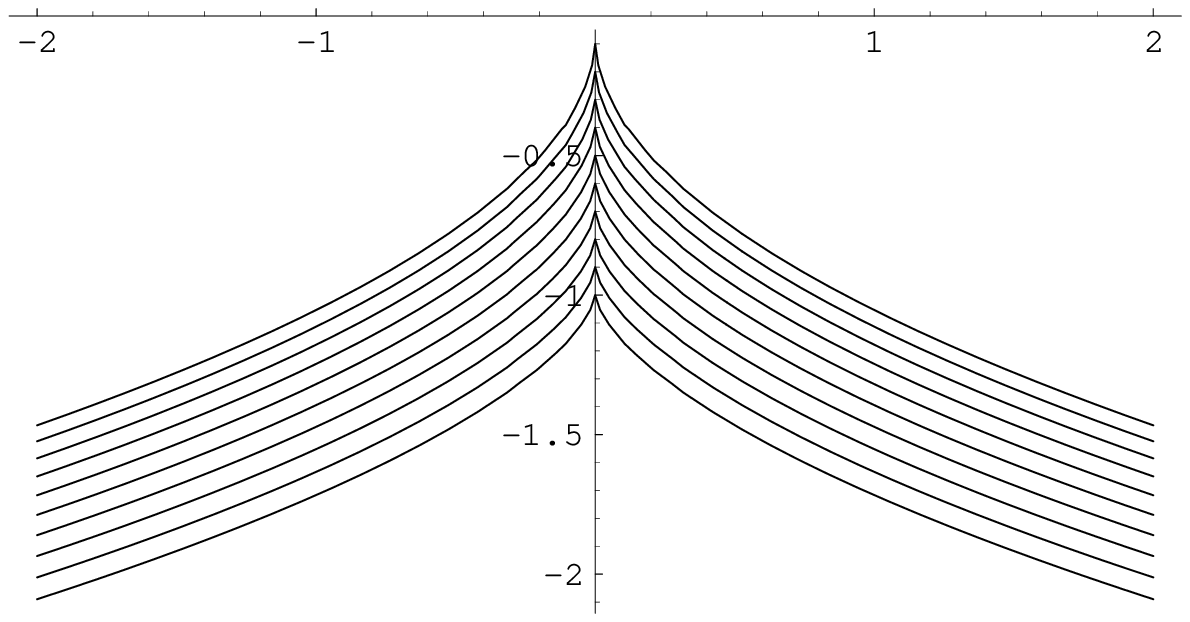}
\quad\quad\quad\quad\quad\quad\quad\quad\quad
\includegraphics[width=28mm]{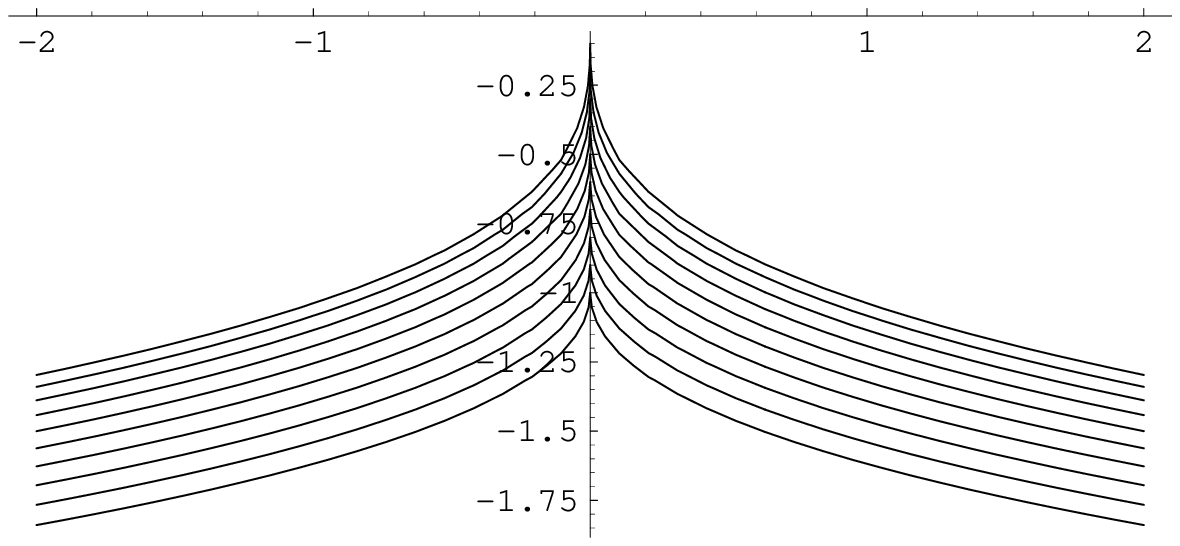}
\end{center}
\caption{Curves generating lower Kummer shapes for  1:-1, 3:-1 and 4:-2 resonance}
\label{fig4}
\end{figure}

The {\bf Kummer shapes} in $n:-m$ resonance are the unbounded surfaces defined by 
the equation (see Figures \ref{fig3} and \ref{fig4}):
\begin{equation*}
x^2+y^2-\left(\frac{z+c}{n}\right)^m\left(\frac{z-c}{m}\right)^n=0,
\quad |z|>c,
\end{equation*}
where $c$ is a positive constant.
Those with $n$ and $m$ of the same parity have two connected components, the others are connected (see Figure \ref{fig7}).
They are obtained by rotating around the $z$ axis the algebraic curve 
$$
y^2=\left(\frac{z+c}{n}\right)^m\left(\frac{z-c}{m}\right)^n,\quad |z|>c.
$$
Let $\Psi\in\F(\RR^4)$ be given by
\begin{equation}\label{psi}
\Psi(x,y,z,r)=x^2+y^2-\left(\frac{z+r}{n}\right)^m\left(\frac{z-r}{m}\right)^n.
\end{equation}
By slicing with hyperplanes $r=c$  of $\RR^4$ that part of the hypersurface $\Psi=0$ included in the union of the halfspaces $z>r$ and $z<-r$, one obtains these Kummer shapes.

\begin{lemma}\label{hardbis}
The unbounded Kummer shapes can be expressed as level sets of a smooth function $C_-$ defined on an open subset of 
$\RR^3\setminus Oz$: 
\begin{equation}\label{p}
B=\{(x,y,z)\in\RR^3\setminus Oz|n^mm^n(x^2+y^2)<z^{n+m}\}.
\end{equation}
\end{lemma}

\begin{proof}
We have to show that there exists a smooth function $C_-:B\to\RR$ which satisfies
\[
\Psi(x,y,z,C_-(x,y,z))=0,\quad C_-(x,y,z)<|z|.
\]
In other words the hypersurface $\Psi=0$ coincides with the graph of the function $C_-$,
on the union of halfspaces $z>r$ and $z<-r$.

To prove this, we observe that $\Psi$ and $B$ are rotationally symmetric in $(x,y)$,
so the problem reduces to proving the existence and uniqueness of a smooth function $c_-$ on $(0,\infty)\x\RR$ such that 
\[
y^2=\left(\frac{z+c_-(y,z)}{n}\right)^m\left(\frac{z-c_-(y,z)}{m}\right)^n,\quad |z|>c_-(y,z),
\]
provided $n^mm^ny^2<z^{n+m}$.
Then $C_-(x,y,z)=c_-(\sqrt{x^2+y^2},z)$ is a smooth function with level sets the unbounded Kummer shapes.

The polynomial function \[ p(r)=\left(\frac{z+r}{n}\right)^m\left(\frac{z-r}{m}\right)^n-y^2, \]
with coefficients depending smoothly on $(y,z)\in(0,\infty)\x\RR$,
has at least one zero in the interval $(0,|z|)$
because $p(0)=\frac{z^{n+m}}{n^mm^n}-y^2>0$ and $p(|z|)=-y^2<0$.
The existence is clear, but for the uniqueness one has to consider separately the two cases $n\ge m$ and $n<m$.
In the first case $p$ is a monotone decreasing function on  $(0,z)$, 
in the second case there is a critical point $r_0\in(0,|z|)$ of $p$, with $p(r_0)>0$,
and $p$ is monotone increasing on $(0,r_0)$ and monotone decreasing on $(r_0,|z|)$. 
In conclusion
there is a unique zero of $p$ in the interval $(0,|z|)$, 
denoted by $c_-(y,z)$.
This gives the smooth function $c_-$ we were seeking.
\end{proof}

The unbounded Kummer shapes in $n:-m$ resonance are symplectic leaves of the Poisson manifold $(B,\pi_\w)$,
where $B$ is defined in \eqref{p} and the bivector field $\pi_\w$ is associated to the vector field
\begin{align}\label{w}
\w&:=\nabla_{(x,y,z)}\Psi\big|_{r=C_-}=\left(2x,2y,-(x^2+y^2)\left(\frac{m}{C_-(x,y,z)+z}-\frac{n}{C_-(x,y,z)-z}\right)\right).
\end{align}
Indeed, $\w=g\nabla C_-$ for $g$ the  nowhere zero function $g=-{\pa_r\Psi}|_{r=C_-}$ on $B$.

As for the $n:m$ resonance, but now with the roles of $Z$ and $R$ switched ($Z_-=R$ and $R_-=Z$), we find that 
\[
X^2+Y^2
=\left(\frac{Z_-+R_-}{n}\right)^m\left(\frac{Z_--R_-}{m}\right)^n
\]
This means that $\Psi\o(X,Y,Z_-,R_-)=0$, so that $C_-\o(X,Y,Z_-)=R_-$ on the open set 
\begin{equation}\label{d}
D=\left\{(a_1,a_2)\in(\CC^2\setminus\{0\})^2: 
(n|a_1|^2)^m(m|a_2|^2)^n<\left(\frac{n}{2}|a_1|^2+\frac{m}{2}|a_2|^2\right)^{n+m}\right\}.
\end{equation}
The inequality defining $D$ comes from $n^mm^n(X(\a)^2+Y(\a)^2)<Z_-(\a)^{n+m}$,
a necessary condition for the existence of $C_-(X(\a),Y(\a),Z_-(\a))$.

\begin{lemma}\label{l2}
For $\{\ ,\ \}_-$ the Poisson bracket on $\CC^2$ induced by the symplectic form $\om_-$, the following identities hold:
\begin{align*}
\{Z_-,X-iY\}_-&=2imn(X-iY)\\
\{X,Y\}_-&=-mn(X^2+Y^2)\left(\frac{m}{R_-+Z_-}-\frac{n}{R_--Z_-}\right).
\end{align*}
\end{lemma}

\begin{proof}
Using the fact that on $\CC$ we have $\{\bar z^n,z^m\}=2imn\bar z^n z^m$, and $\{|z|^2,z^n\}=2inz^n$ 
as well as $\{|z|^2,\bar z^n\}=-2inz^n$,
we compute
\begin{align*}
\{Z_-,X-iY\}_-&=\frac n2\bar a_2^n\{|a_1|^2,a_1^m\}_- 
+\frac m2 a_1^m \{|a_2|^2\bar a_2^n\}_- \\
&=\frac n2(2im)a_1^m\bar a_2^n-\frac m2(-2in)a_1^m\bar a_2^n=2imn(X-iY)
\end{align*}
and
\begin{align*}
\{X,Y\}_-&=\frac{i}{2}\{\bar a_1^ma_2^n,a_1^m\bar a_2^n\}_-
=\frac{i}{2}\{\bar a_1^m,a_1^m\}_-|a_2|^{2n}-\frac{i}{2}\{\bar a_2^n,a_2^n\}_-|a_1|^{2m}\\
&=\frac{i}{2}(2im^2|a_1|^{2m-2})|a_2|^{2n}-\frac{i}{2}(-2in^2|a_2|^{2n-2})|a_1|^{2m}\\
&=-|a_1|^{2m}|a_2|^{2n}\left(\frac{m^2}{|a_1|^2}+\frac{n^2}{|a_2|^2}\right)
=-mn(X^2+Y^2)\left(\frac{m}{R_-+Z_-}-\frac{n}{R_--Z_-}\right).
\end{align*}
\end{proof}

\begin{proposition}\label{p3}
The map $\Pi_-=(X,Y,Z_-):D\subset\CC^2\to B\subset\RR^3$ is a Poisson map with respect to 
the symplectic form $\omega_-$ on $\CC^2$ and
the Poisson bivector field $\pi_{mn\w}$ on $B$.
\end{proposition}

\begin{proof}
From lemma \ref{l2} we have that:
\begin{align*}
\{Y,Z_-\}_-&=2mnX\\
\{Z_-,X\}_-&=2mnY\\
\{X,Y\}_-&=-mn(X^2+Y^2)\left(\frac{m}{R_-+Z_-}-\frac{n}{R_--Z_-}\right)
\end{align*}
Knowing that
\[
\pi_{mn\w}=2mnx\pa_y\wedge\pa_z+2mny\pa_z\wedge\pa_x-mn(x^2+y^2)\left(\frac{m}{C_-(x,y,z)+z}
-\frac{n}{C_(x,y,z)--z}\right)\pa_x\wedge\pa_y,
\]
the result follows from the functional identity $C_-\o(X,Y,Z_-)=R_-$.
\end{proof}

\begin{theorem}\label{t2}
The pair of momentum maps
\[
{\RR\stackrel{R_-}{\longleftarrow}(D,\om_-)
\stackrel{\Pi_-}{\longrightarrow}(B,\pi_{mn\w})}
\] 
is a dual pair for all pairs $(m,n)$ of nonzero natural numbers, with $B$ and $D$ given in \eqref{p} and \eqref{d}.
\end{theorem}

\begin{proof}
We know already that both $R_-$ and $\Pi_-$ are Poisson maps.
We have to show the dual pair property 
\begin{equation}\label{dpdp}
\ker TR_-=(\ker T\Pi_-)^{\om_-}.
\end{equation}
The proof is similar to that of Theorem \ref{t1}.
The symplectic orthogonal for $\om_-$ and the Riemannian orthogonal for the canonical Riemannian metric
$g$ on $\CC^2$
are related by: $(a_1,a_2)^\perp=(ia_1,-ia_2)^{\om_-}$ because of the identity
$\om_-(\a,\b)=g((a_1,a_2),(ib_1,-ib_2))$
for all $\a=(a_1,a_2)$ and $\b=(b_1,b_2)$ on $\CC^2$.
Thus we get 
\begin{equation}\label{two}
\ker T_{\a}R_-=(na_1,-ma_2)^\perp=(nia_1,mia_2)^{\om_-}.
\end{equation}
As in the proof of Theorem \ref{t1} one sees that 
\[
(nia_1,mia_2)\in\ker T_{\a}\Pi_-=\ker T_{\a}X\cap\ker T_{\a}Y\cap\ker T_{\a}Z_-
\]
The kernel of $T_{\a}\Pi_-$ being 1-dimensional, it must be generated by $(nia_1,mia_2)$.
We get
\[
(\ker T_{\a}\Pi_-)^{\om_-}=(nia_1,mia_2)^{\om_-},
\]
which, together with \eqref{two}, ensures the dual pair property \eqref{dpdp}.
\end{proof}

The symplectic leaf correspondence theorem for dual pairs, applied to the $n:-m$ resonance, says
that for each $c\in\RR$, the symplectic leaf $\{c\}$ of $\RR$ corresponds to the symplectic leaf
$\Pi(R^{-1}(c))$ of $B$, \ie to the unbounded Kummer shape $C_-(x,y,z)=c$, because $C_-\o\Pi_-=R_-$.

\section{Conclusions}

A Hamiltonian system having a number of independent integrals of motion bigger than 
the dimension of its phase space is called superintegrable.
More precisely we are given a symplectic $2d$-dimensional symplectic manifold $(M,\om)$,
and a submersion $f=(f_1,\dots,f_{2d-n}):M\to \RR^{2d-n}$ with compact connected fibers,
with two properties:
\begin{enumerate}
\item $\{f_i,f_j\}=\pi_{ij}\o f$ for $\pi_{ij}:B\subset\RR^{2d-n}\to\RR$, $i,j=1,\dots,2d-n$,
\item $\rank (\pi_{ij})=2d-2n$.
\end{enumerate}
The Michenko-Fomenko theorem \cite{MF} states that under these circumstances 
the fibers of $f$ are diffeomorphic to the $n$-dimensional torus,
and locally there exist generalized action-angle coordinates
$(p,q,a,\al)$ on the symplectic manifold $M$ ($p$ and $q$ have $d-n$ components, while the actions $a$ and the angles $\al$ have $n$ components), that is $\om=dp\wedge dq+da\wedge d\al$.
Integrable systems with $d$ integrals of motion in involution are obtained for $n=d$.

As explained in \cite{Fasso}\cite{OrRa04}, the two conditions have a geometric interpretation.
The functions $\pi_{ij}$ are the components of a Poisson bivector field $\pi_B$ on the open subset 
$B\subset\RR^{2d-n}$, such that
the map $f:(M,\om)\to(B,\pi_B)$ is Poisson.
By dimension counting follows that $\rank\pi_B=\dim M-\dim\ker Tf-\dim(\ker Tf\cap(\ker Tf)^\om)$,
so condition 2. implies that $\ker Tf\subset(\ker Tf)^\om$, which means that the submersion $f$
has isotropic fibers (the isotropic tori).

Condition 1. also ensures that the orthogonal distribution $(\ker Tf)^\om$ is integrable.
This follows from its involutivity: for the local basis $X_{f_1},\dots,X_{f_{2d-n}}$ of $(\ker Tf)^\om$
consisting of the Hamiltonian vector fields with Hamiltonian functions given by the $2d-n$  integrals of motion,
all commutators $[X_{f_i},X_{f_j}]=X_{\{f_i,f_j\}}=X_{\pi_{ij}\o f}$
are again sections of $(\ker Tf)^\om$.
The other way around, the obvious integrability of $\ker Tf$ ensures that
there is a Poisson bivector field $\pi_A$ on the space $A$ of leaves of 
the integrable distribution $(\ker Tf)^\om$ ($A$ is assumed to be a manifold)
such that the projection on the space of leaves $p:(M,\om)\to(A,\pi_A)$ is Poisson.

The dual pair of Poisson maps associated to the superintegrable system is
$$
(A,\pi_A)\stackrel{p}{\longleftarrow}(M,\om)\stackrel{f}{\longrightarrow}(B\subset\RR^{2d-n},\pi_B).
$$
The angles $\al$ are coordinates on the fibers of $f$, while the actions $a$ are local coordinates on $A$
and $(p,q)$ are local coordinantes on the symplectic leaves of $B\subset\RR^{2d-n}$.
These two Poisson maps coincide in the integrable case $n=d$.

The dual pair for the superintegrable system of two uncoupled oscillators in $m:n$ resonance 
is the one presented in Theorem \ref{t1} for positive $n$, resp. Theorem \ref{t2} for negative $n$.
The functions $X,Y,Z$ \eqref{xyz}, resp. $X,Y,Z_-$ \eqref{nega}, are the three independent integrals of motion
on the 4-dimensional symplectic manifold $(M=(\CC\setminus\{0\})^2,\om)$, resp. $(M=\{(a_1,a_2)\in(\CC\setminus\{0\})^2: (n|a_1|^2)^m(m|a_2|^2)^n<\left(\frac{n|a_1|^2+m|a_2|^2}{2}\right)^{n+m}\},\om_-)$, where $\om$ is the opposite of the canonical symplectic form on $\CC^2$,
and $\om_-=-\frac i2(da_1\wedge d\bar a_1-da_2\wedge d\bar a_2)$.

The Poisson manifold $B$ is an open subset of $\RR^3$: $B=\RR^3\setminus Oz$
with Poisson bivector field 
$\pi_B=2mnx\pa_y\wedge\pa_z+2mny\pa_z\wedge\pa_x-mn(x^2+y^2)\left(\frac{m}{C(x,y,z)+z}-\frac{n}{C(x,y,z)-z}\right)\pa_x\wedge\pa_y$,
where $C$ is a smooth function on $B$ implicitly defined by $x^2+y^2-\left(\frac{C(x,y,z)+z}{n}\right)^m\left(\frac{C(x,y,z)-z}{m}\right)^n=0$ and
$|z|<C(x,y,z)$,
resp. $B=\{(x,y,z)\in\RR^3\setminus Oz|n^mm^n(x^2+y^2)<z^{n+m}\}$ with Poisson bivector field $\pi_B=2mnx\pa_y\wedge\pa_z+2mny\pa_z\wedge\pa_x-mn(x^2+y^2)\left(\frac{m}{C_-+z}-\frac{n}{C_--z}\right)\pa_x\wedge\pa_y$, where $C_-$ is a smooth function on $B$ implicitly defined by
$x^2+y^2-\left(\frac{z+C_-(x,y,z)}{n}\right)^m\left(\frac{z-C_-(x,y,z)}{m}\right)^n=0$ and
$C_-(x,y,z)<|z|$.

The symplectic leaves of the Poisson manifold $(B,\pi_B)$ are the Kummer surfaces,
bounded for positive $n$, resp. unbounded for negative $n$.


{\footnotesize

\bibliographystyle{new}

\addcontentsline{toc}{section}{References}
\begin{thebibliography}{300}

\bibitem[Cushman and Rod(1982)]{CR}
R. Cushman and D.~L. Rod [1982],
Reduction of the semisimple 1:1 resonance,
\textit{Physica D} \textbf{6}, 105--112.

\bibitem[Elipe(2000)]{Elipe}
A. Elipe [2000],
Complete reduction of oscillators in resonance $p:q$,
\textit{Phys. Rev. E} \textbf{61}, 6477--6484.

\bibitem[Fass\`{o}(2005)]{Fasso}
F. Fass\`{o} [2005],
Superintegrable Hamiltonian systems: geometry and perturbations,
\textit{Acta Appl. Math.} \textbf{87}, 93--121.

\bibitem[Gay-Balmaz and Vizman(2009)]{GBV10}
\textit{F. Gay-Balmaz and C. Vizman}, Dual pairs in fluid dynamics,
preprint 2009.

\bibitem[Holm(2008)]{Holm}
D.~D. Holm [2008],
\textit{Geometric Mechanics Part 1: Dynamics and Symmetry}, World Scientific, London. 


\bibitem[Iwai(1985)]{Iwai}
T. Iwai [1985],
On reduction of two degrees of freedom Hamiltonian systems by an $S^1$ action and $SO_0(1,2)$
as a dynamical group,
\textit{J. Math. Phys.} \textbf{26}, 885--893.


     \bibitem[Kummer(1976)]{Ku1976}
     Kummer, M. [1976],
     On resonant nonlinearly coupled
oscillators with two equal frequencies.
     {\it  Commun. Math. Phys.} {\bf 48}, 53-79 (1976).

     \bibitem[Kummer(1978)]{Ku1978}
     Kummer, M. [1978],
On resonant classical Hamiltonians
with two equal frequencies.
     {\it  Commun. Math. Phys.} {\bf 58}, 85-112.

     \bibitem[Kummer(1981)]{Ku1981}
     Kummer, M. [1981],
     On the construction of the reduced phase space of a Hamiltonian system with symmetry. {\it Indiana Univ. Math. J.} {\bf 30}, 281-291.
     
     \bibitem[Kummer(1986)]{Ku1986}
Kummer,  M.  [1986] in {\it Local and Global Methods in Nonlinear
Dynamics}, Lecture Notes in Physics Vol. {\bf252}, edited by 
A. V. S\'aenz, Springer-Verlag,  New York, pp. 19-31.

\bibitem[Marsden(1987)]{Ma}
J.E. Marsden [1987], 
Generic Bifurcation of Hamiltonian Systems with Symmetry, \textit{appendix to Golubitsky and Stewart, Physica D}, \textbf{24}, 391--405.

\bibitem[Mishenko and Fomenko(1978)]{MF}
Mishenko, A.~S. and Fomenko, A.~T. [1978],
Generalized Liouville method of integration of Hamiltonian systems, Funct. Anal. Appl. \textbf{12}, 113--121.

\bibitem[Ortega and Ratiu(2004)]{OrRa04}
Ortega, J.-P. and T.~S. Ratiu [2004], \textit{Momentum maps and
Hamiltonian reduction}, Progress in Mathematics (Boston, Mass.)
\textbf{222} Boston,  Birkh\"auser.

\bibitem[Weinstein(1983)]{We83}
Weinstein, A. [1983], The local structure of Poisson manifolds, \textit{J. Diff. Geom.} \textbf{18}, 523--557.


\end{thebibliography}

}

\end{document}